\documentclass[12pt,a4paper]{article}
 
 \usepackage{amssymb, graphics, setspace}
\usepackage{romanbar}
\usepackage{bm}
\usepackage{subfig}
\newcommand{\mathsym}[1]{{}}
\newcommand{\unicode}[1]{{}}

\usepackage[pdftex]{hyperref}
\usepackage{amssymb}
\usepackage[centertags]{amsmath}

\usepackage{amsthm}
\usepackage{graphicx}
\graphicspath{ {./Pictures/} }

\usepackage{xcolor}

\newtheorem{thm}{Theorem}

\newtheorem{pro}{Proposition}
\newtheorem*{pro*}{Proposition}

\newtheorem{rem}{Remark}
\newtheorem*{rem*}{Remark}
\newtheorem{cor}{Corollary}

\renewcommand{\proof}{\textsc{proof. }}
\newcommand{\eproof}{\hfill $\Box$}
\newcommand{\bal}{{\bm \alpha}}
\newcommand{\bbe}{{\bm \beta}}
\newcommand{\bom}{{\bm \omega}}
\newcommand{\bphi}{{\bm \varphi}}
\newcommand{\bpsi}{{\bm \psi}}
\newcommand{\bI}{\mathbf I}
\newcommand{\bJ}{\mathbf J}

\makeindex

\title{On bandgap sensitivity to three-to-one 
internal resonances
 between acoustic and optical waves in metamaterials}

\author{L. Di Gregorio$^{*}$ and W. Lacarbonara$^{*}$\\
\small{laura.digregorio@uniroma1.it, \, walter.lacarbonara@uniroma1.it}\\
\footnotesize{$^{*}$Department of Structural and Geotechnical
Engineering, Sapienza University of Rome,}\\
\footnotesize{Via Eudossiana 18, 00184, Rome,
Italy}}

\begin{document}

\allowdisplaybreaks
\maketitle

{\abstract{We investigate the nonlinear equations
governing wave propagation
 across a metamaterial consisting of a cellular periodic  structure hosting resonators with linear and cubic springs. The resulting system of two coupled equations with cubic nonlinearity is Hamiltonian, with the origin being an elliptic equilibrium characterized by two distinct linear frequencies associated with the acoustic and optical wave modes. Understanding this complex wave propagation problem requires the explicit derivation of analytical formulas for the nonlinear frequencies as functions of the relevant phy\-si\-cal parameters. In the small wave amplitude regime and away from internal resonances, by using Hamiltonian Perturbation Theory, we obtain the first-order nonlinear correction to the linear frequencies. 
 Additionally, we analytically estimate the threshold for the
applicability of the perturbative solution, identifying the maximum
admissible amplitude at which the 
obtained formulas remain valid.
  This threshold decreases significantly in the presence of internal resonances, particularly when the ratio between the optical and acoustic frequencies is close to 3
  due to a 
 3:1 resonance which affects  the correction. Furthermore, we analytically evaluate the remainder after the first step of the perturbation scheme. The specific application deals with the wave propagation equations obtained for metamaterial honeycombs with periodically distributed nonlinear one-dof resonators. Our analysis provides quantitative insights into how internal resonances affect the bandgap. Specifically, we identify a 3:1 resonant curve in the parameter space represented by the resonator modal mass and stiffness. We also identify a region   away from the
 above resonant curve where the nonlinear bandgap for softening resonators springs is significantly larger than the linear counterpart.

 }}

 \bigskip
\noindent
{\footnotesize{\textbf{Acknowledgments} Project ECS 0000024 Rome Technopole, CUP  B83C22002820006, National Recovery and Resilience Plan (NRRP) Mission 4 Component 2 Investment 1.5, funded by the European Union - NextGenerationEU.}}

\bigskip

\noindent
{\footnotesize{\textbf{Funder} Project funded under the National Recovery and Resilience Plan (NRRP), Mission 4 Component 2 Investment 1.5 - Call for tender No. 3277 of 30 December 2021 of the Italian Ministry of University and Research funded by the European Union - NextGenerationEU.}}
\tableofcontents

\section{Introduction}

This study was inspired by \cite{SW23jsv} and \cite{SW23mssp}, where the goal was to utilize the resonators nonlinearity to create and optimize targeted flexural wave bandgaps in 2D locally resonant periodic structures. Therein, the authors provided design conditions for the softening and hardening resonators nonlinearity in metamaterial honeycombs, using asymptotic solutions.

Regarding the model, earlier investigations into 2D locally resonant periodic structures with linear resonators can be found in \cite{Liu21, Fan21, M22, Guo22, B16, Comi18, Wang21}. Various lattice structures and resonator architectures in different shapes have been proposed to control wave propagation in \cite{Cai22}. Two-dimensional weakly nonlinear lattices were asymptotically investigated in \cite{L19jsv}, with a special focus on invariant waveforms and stability analysis. Using the method of multiple scales, Bukhari and Barry \cite{B20} examined a nonlinear metamaterial consisting of a nonlinear chain with multiple nonlinear local resonators, revealing the relationship between the space-time domain and nonlinear dispersion properties.
A key role in studies of 2D cellular structures hosting resonators is played by the homogenization approach, which is implemented to obtain equivalent continuum plate-like models, as in \cite{G97}, \cite{S18}.
The effects of dissipation in a 1D hosting medium made of linear springs with a periodic distribution of nonlinear, hysteretically damped resonators were previously addressed using an extended Hamiltonian approach in \cite{F22}.

In terms of bandgap analysis, in the context of linear behaviors, the work in \cite{M19} investigates the dispersion properties of cellular locally resonant metamaterials with various geometries in both two and three dimensions. This study employs finite element modal analysis on the unit cell with Bloch-Floquet boundary conditions. Similarly, \cite{J19} explores the bandgaps of locally resonant materials with different plate-like resonators using a finite element model based on plate theory to analyze the dispersion relations.

\noindent
For the nonlinear case, the study in \cite{Wi19} reports a global Sobol sensitivity analysis using Monte Carlo integration to identify the most important parameters in designing frequency band gaps of periodic materials. The work in \cite{My22} introduces a new nonlinear metamaterial featuring a tunable bandgap at extremely low frequencies, starting from zero. This metamaterial achieves nonlinearity through geometric changes in the bounding springs connecting each mass to fixed foundations, offering dual tunability through wave amplitude and spring length variation. Additionally, a detailed perturbation analysis of wave dispersion in one-dimensional, discrete, nonlinear periodic structures is discussed in \cite{Na20}. 
The present work aims to investigate  analytically the  sensitivity of the {\sl fundamental  bandgap} with respect to the system parameters and the role played by three-to-one interactions between the acoustic and optical waves surrounding the bandgap.

\subsubsection*{The model}

\begin{figure}
\center
\includegraphics[width=6cm,angle=0.7,keepaspectratio]{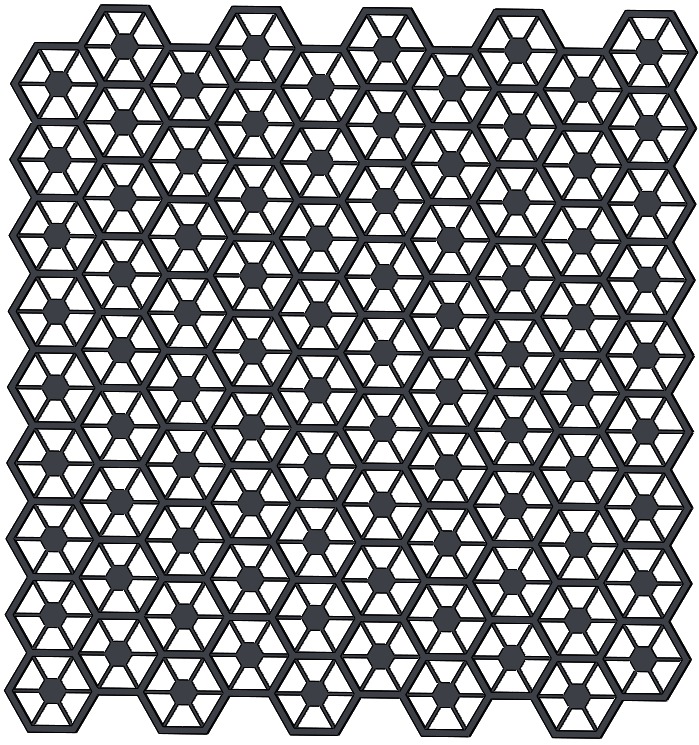}
\caption{Schematic view of a 2D metamaterial with periodically 
distributed resonators, proposed in \cite{SW23jsv}.}
\label{Resonators_model}
\end{figure}
Let us briefly illustrate the model introduced in 
\cite{SW23jsv}.
Figure \ref{Resonators_model}
 shows a schematic view of the 2D metamaterial with the periodically distributed spider-web-inspired resonators.
 The 2D hosting honeycomb was
  modeled as an equivalent orthotropic plate. 
 Each multi-frequency resonator is equivalent 
 to a multi-mass-spring system resulting from the modal reduction of the infinite-dimensional resonator structure. In this work, the 
 resonators with the
  central mass will be represented by a single mass-spring system instead of a set of mass-spring systems. The modal reduction is performed via the Galerkin projection method employing the fundamental mode shape of the resonator.

\begin{figure}[h!]
			\centering
\includegraphics[width=5cm,height=5cm,keepaspectratio]{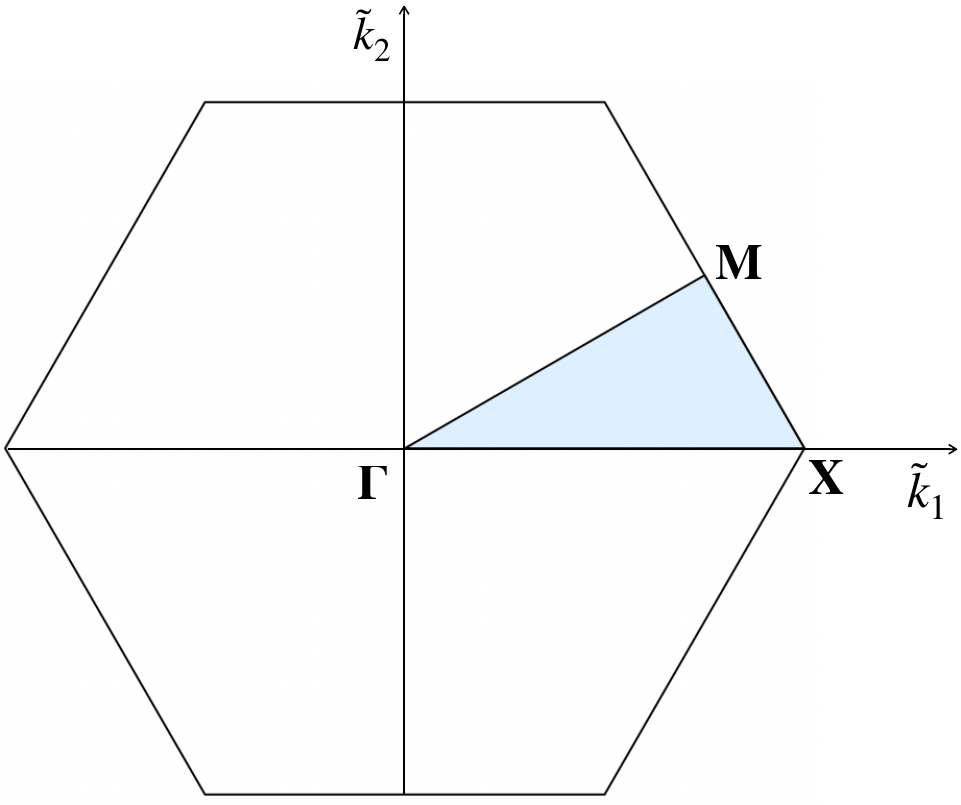}
\caption{The irreducible Brillouin triangle 
$\triangle:=\bf{\Gamma}\overset{{\triangle}}{\bf{X}} \bf{M}$ with $\mathbf{\Gamma}=(0,0)$, $\mathbf{X}=(\frac{4}{3}\pi,0)$, $\mathbf{M}=(\pi,\frac{\pi}{\sqrt{3}})$.}
\label{Brillouin}
		\end{figure}

The adopted plate theory  (see \cite{W}) with the elastic constants of the equivalent, homogenized  orthotropic material  describes the motion of the honeycomb with the attached resonators.
By  the Floquet-Bloch Theorem, which states that the solutions of the {corresponding linear} periodic  resonators-plate system are quasi-periodic in space with the fundamental periodicity provided by the lattice period,  
the plate equation of motion can be projected onto the unit cell domain (i.e., the periodically repeated lattice unit).
Then, one obtains the following
 system of second order ODEs:
\begin{equation}\label{firenze}
 \left(
\begin{array}{cc}
	\tilde{M}_H(\tilde{k}_1,\tilde{k}_2)&\tilde{M}\\
       \tilde{M}& \tilde{M}
	\end{array}	
\right)
\left(
		\begin{array}{c}
			\ddot{\tilde w}_0 \\
			\ddot{\tilde z}_0
		\end{array}
		\right)
		+
		 \left(
\begin{array}{cc}
	\tilde{K}_H(\tilde{k}_1,\tilde{k}_2)&0\\
       0& \tilde{K}
	\end{array}	
\right)
		\left(
		\begin{array}{c}
			{\tilde w}_0 \\
			{\tilde z}_0
		\end{array}
		\right)
		=
		-\left(
		\begin{array}{c}
			 0 \\
			\tilde N^{(3)} {\tilde z}_0^3
		\end{array}
		\right)\,,
\end{equation}
where
${\tilde w}_0$ and 
			${\tilde z}_0$
denote the nondimensional plate deflection and  resonator relative motion at the origin of the fixed frame; $\tilde M$, $\tilde K$ are the modal mass and stiffness of the resonator and
	\begin{equation}\label{coefficients3}
		\tilde{M}_H (\tilde{k}_1,\tilde{k}_2) 
		:=
		 \frac{4 \sqrt{3} \sin \left(\frac{\tilde{k}_1}{2}\right) \sin \left(\frac{1}{4}  \left(\tilde{k}_1+\sqrt{3} \tilde{k}_2\right)\right)}{\tilde{k}_1 \left(\tilde{k}_1+\sqrt{3} \tilde{k}_2\right)} \,,
	\end{equation}
\begin{equation}\label{coefficients2}
	\begin{split}
		&
\tilde{K}_H(\tilde{k}_1,\tilde{k}_2)= 
\tilde{K}_H(\tilde{k}_1,\tilde{k}_2;
\tilde{D}_{12},\tilde{D}_{66},\tilde{D}_{22})
:=\tilde{M}_H(\tilde{k}_1,\tilde{k}_2) \left[\tilde{k}_1^4  +2 \tilde{k}_1^2 \tilde{k}_2^2 (\tilde{D}_{12} + 2\tilde{D}_{66}) + \tilde{k}_2^4 \tilde{D}_{22} \right]
	\end{split}
\end{equation}
are the nondimensional mass and stiffness
expressed as functions of the nondimensional
wave numbers $(\tilde{k}_1,\tilde{k}_2)$,
which vary in the range provided by the irreducible Brillouin zone (see Figure \ref{Brillouin}).
Moreover,	
$$
\tilde D_{12}=0.0815599,\qquad \tilde D_{22}=12.48, 
\qquad\tilde D_{66}=0.0000247357\,,
$$
 are the nondimensional bending stiffness coefficients of the orthotropic plate model; finally
 $\tilde{N}^{(3)}$ is the nondimensional 
spring constant.

\medskip
Actually, for sake of generality, we  consider a more
 general system of ODEs; that is
\begin{equation}\label{autostrada}
\mathtt M\left(
		\begin{array}{c}
			\ddot v \\
			\ddot y
		\end{array}
		\right)
		+\mathtt K\left(
		\begin{array}{c}
			v \\
			y
		\end{array}
		\right)
		=
		-\left(
		\begin{array}{c}
			 M_3v^3 \\
			N_3 y^3
		\end{array}
		\right)\,,
\end{equation}
where $v(t)$ and $y(t)$
are unknown scalar functions,
$M_3$ and $N_3$ are real coefficients,
$\mathtt M$ is a symmetric 
positive-definite $2\times 2 $ real mass matrix
and $\mathtt K$ is a diagonal 
positive-definite $2\times 2 $ real stiffness matrix.


Note that 
 \eqref{firenze} is a particular
 case of \eqref{autostrada} taking
 $v=\tilde{w}_0$, $z=\tilde{z}_{0}$,
 $M_3=0,$ $N_3=\tilde{N}^{(3)}$
 and
 \begin{equation}\label{frittata}
 \mathtt M=
 \left(
\begin{array}{cc}
	\tilde{M}_H(\tilde{k}_1,\tilde{k}_2)&\tilde{M}\\
       \tilde{M}& \tilde{M}
	\end{array}	
\right)\,,\qquad
 \mathtt K=
 \left(
\begin{array}{cc}
	\tilde{K}_H(\tilde{k}_1,\tilde{k}_2)&0\\
       0& \tilde{K}
	\end{array}	
\right)\,,
 \end{equation}
 with $\tilde{M}_H(\tilde{k}_1,\tilde{k}_2)$
 and 
 $\tilde{K}_H(\tilde{k}_1,\tilde{k}_2)$
 being defined in 
 \eqref{coefficients3} and \eqref{coefficients2},
 respectively.

\medskip

\subsubsection*{Main results}

We are interested in small
amplitude solutions
of \eqref{autostrada}.  Thus, in the first
approximation, the system is linear
with frequencies 
$\omega_-$ and $\omega_+$
while the nonlinearity is a third order
perturbation.
If the linear frequencies 
are non vanishing, {\sl away from the 3:1 resonance condition
$3\omega_-= \omega_+$},
the system can be integrated,
e.g. by the multiple scales approach,
up to a smaller nonlinear remainder
 of fifth order, see \cite{SW23jsv}.
In particular, in \cite{SW23jsv} the nonlinear frequencies
of the truncated system (obtained
disregarding the fifth order perturbation) were 
explicitly given as functions of the
initial amplitudes, see formula
\eqref{formulaSL1}. Moreover,
the effects of the nonlinear local resonances on the bandgap
were thoroughly explored.
\\
It turns out that nonlinearity can significantly affect the bandgap,
either expanding or reducing it, depending on the type of nonlinearity
and the wave amplitudes. In this study, we examine the impact of
resonances on this process.
According to our analysis, it follows that,
in order to maximise the bandgap gain triggered by the nonlinear local resonance,
one should choose
the parameters in 
such a way that 
the quantity 
$|3 \omega_--\omega_+|$
is large.
 Indeed,
 this allows to enlarge
the size of the admissible initial
amplitudes, namely,
 the amplitudes for which the 
perturbative procedure works and formula
\eqref{formulaSL1} makes sense. 
As a consequence, since the increment
of the bandgap
is proportional to the square
of the amplitudes,  
one can
  optimise  the 
boosting effects of the
nonlinearity on the bandgap width.
\\
Finally,
 we analytically evaluate the
 remainder after the first step of the perturbative procedure.
 This enables us to estimate the distance between the third-order
truncated system's solutions and the true solutions of the full
problem over long times.

\medskip

\begin{figure}[h!]
			\centering			
			\includegraphics[width=7cm,height=7
		cm,keepaspectratio]{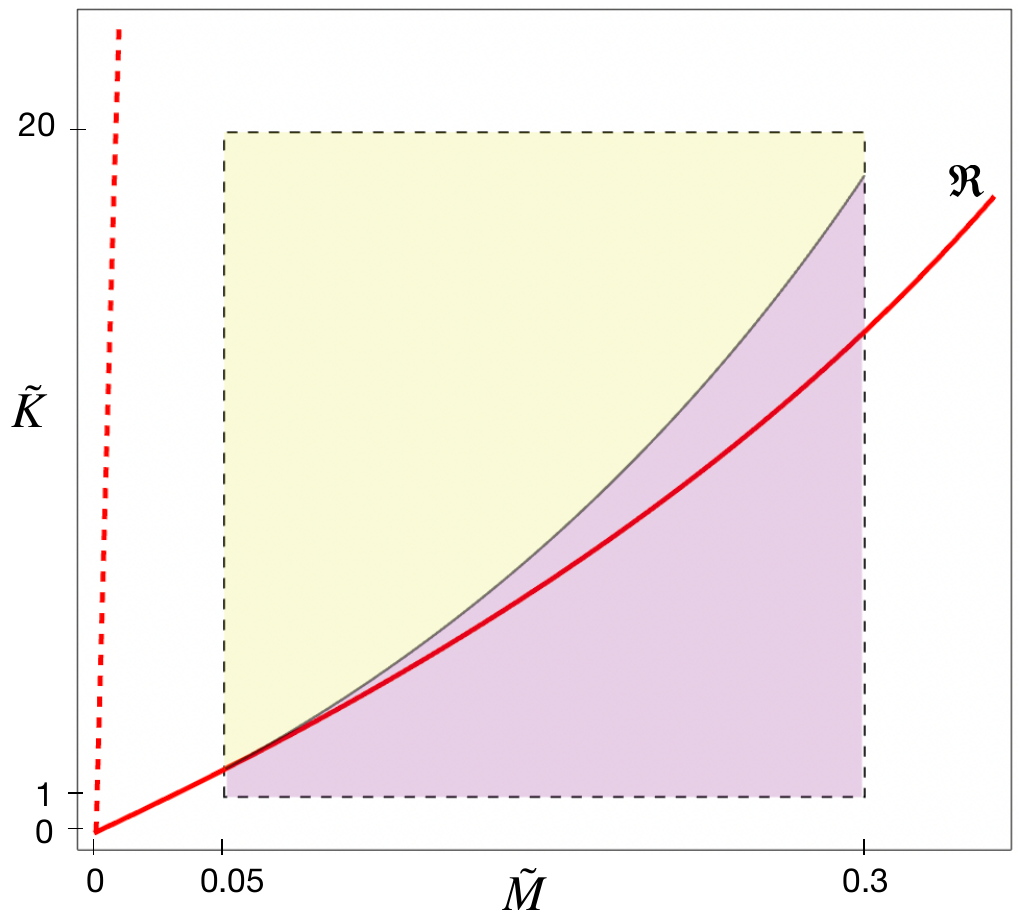}
			\caption{The solid red curve 
			$\mathfrak R$ and the dashed red curve are the 
			 two ${\bf X}$-resonant
			curves in the $(\tilde M,\tilde K)$-plane. The light yellow region within the  reference rectangle $[0.05, 0.3]\times[1, 20]$
			indicates 
			 the ``good'' nonresonant set where the  representation formula \eqref{formulaSL1}
			can be applied.
	In contrast, the light purple region
	shows the ``bad'' set of resonant parameters,
	where a ``resonant'' representation 
	formula is necessary.	
			 }
			 \label{buoni}
		\end{figure}

As an application of our analysis we consider the honeycomb metamaterial
discussed above and carry out the evaluation of its nonlinear bandgap.
A relevant aspect is that the linear frequencies and, particularly the
quantity $|3 \omega_--\omega_+|$, strongly depend, not only 
on the modal parameters $\tilde M$
and  $\tilde K$, but also on the
wave numbers $\tilde k_1$ and $\tilde k_2$ as they vary along the boundary of the
first  IBZ (irreducible Brillouin zone); see Figure \ref{Brillouin}. 
The bandgap is defined as the gap between the minimum of the optical frequency and the maximum of the acoustic
frequency over the IBZ  range. Actually, these extrema are reached on the boundary of the triangle. 
In the linear case, these extremal values are
attained at the vertices
${\bf X}$ and ${\bf\Gamma}$, respectively.
Notably, the vertex 
${\bf X}$ plays a more critical role than ${\bf \Gamma}$.
Indeed, the shift in the maximum
of the acoustic frequency due to the nonlinearity is more significant than the shift 
in the minimum of the optical frequency.
Consequently, we identify two {\sl ${\bf X}$-resonant} curves in the
 $(\tilde M,\tilde K)$-plane, which are curves formed by pairs
 $(\tilde M,\tilde K)$ such that, when
 $(\tilde k_1,\tilde k_2)={\bf X}$, they satisfy the exact 3:1 resonance
 $3\omega_-=\omega_+$ (see the  solid and dashed red curves
 in Figure \ref{buoni}). Since our analysis is restricted 
 to the reference rectangle  $[0.05, 0.3]\times[1, 20]$, we 
 consider only the solid ${\bf X}$-resonant curve shown in the figure.

Next, we identify the set of
 ``good'' pairs $(\tilde M,\tilde K)$ 
within  the rectangle $[0.05, 0.3]\times[1, 20]$ (represented by the light yellow region in Figure
\ref{buoni}), for which
 the maximum/minimum of the nonlinear acoustic/optical
 frequencies on the boundary
of the Brillouin triangle,
 are attained at nonresonant
wave numbers $(\tilde k_1,\tilde k_2)$, i.e.,
at points where the quantity  $|3 \omega_--\omega_+|$
is not small.

Formula \eqref{formulaSL1} is valid in this  ``good'' set, 
allowing us to directly evaluate the bandgap.
In contrast, within the complementary light purple zones in 
Figure
\ref{buoni}, formula \eqref{formulaSL1} is not applicable.
To evaluate the nonlinear bandgap in the light purple
 zones, a new ``resonant'' representation
 formula for the nonlinear frequencies must be used, which is more complex
 and  is derived  in \cite{DL2}.
 
\begin{figure}[h!]
			\centering			
			\includegraphics[width=7cm,height=7cm,keepaspectratio]{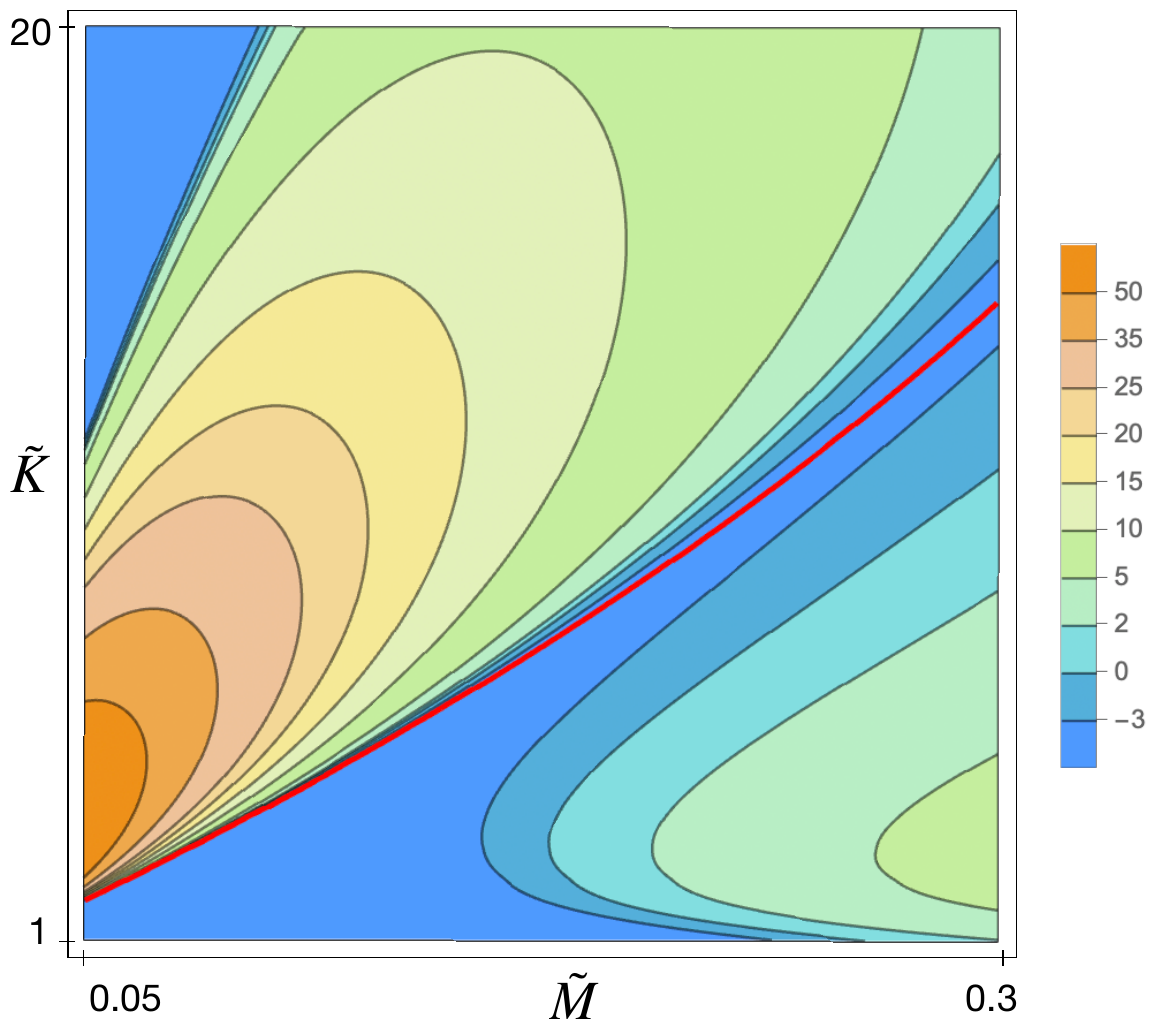} 
			\caption{
		Contour plot of the maximum percentage 		difference between nonlinear and linear bandgap
in the $(\tilde M,\tilde K)$ plane.
The ${\bf X}$-resonant curve  $\mathfrak R$
is plotted in red.
			Here,  $N_3=- 10^4$ (softening spring).
			In this softening case,
			the majority of parameter pairs above $\mathfrak R$
			 result in
			an increase in the bandgap width,
			 while those below $\mathfrak R$ either show a decrease or, at most, a very slight increase.
			 The region where the increase is most pronounced is contained in the set of ``good'' parameters highlighted in light yellow in	Figure \ref{buoni}.}
			 \label{ritz}
					\end{figure}

The final result of our analysis is shown in Figure \ref{ritz}, where the
maximum percentage increase in the nonlinear bandgap, compared to the
linear bandgap, is plotted 
 as the pair $(\tilde M,\tilde K)$
 varies over the rectangle $[0.05, 0.3]\times[1, 20]$
 in the softening case (i.e. the spring constant was set to $N_3=-10^4$).
It is important to note that, in this paper, we 
demonstrate how to derive Figure \ref{ritz} using \eqref{formulaSL1}, only for 
 the pairs $(\tilde M,\tilde K)$ belonging
 to the light yellow set in Figure \ref{buoni}. As mentioned above, for the pairs $(\tilde M,\tilde K)$ belonging
 to the light purple set, we refer to \cite{DL2}.

\textsl{
We emphasize that our analysis is not confined to honeycomb metamaterials alone but is directly applicable to a broad range of problems modeled by two  harmonic oscillators coupled with cubic nonlinearity, as described 
in Eq. \eqref{autostrada}}. Studying  nonlinear oscillations of these coupled harmonic oscillators can offer deeper insights into the dynamic behavior of more complex systems. Furthermore, it could contribute to the development of more effective vibration and wave control strategies. 
\medskip

\subsubsection*{Methodology}

Since the system is conservative,
 we study it as a Hamiltonian 
 system. The origin is an elliptic
 equilibrium and we put the 
 system in {\sl Birkhoff Normal Form} (BNF)
 up to order 4 for the system Hamiltonian (3 in the equations
 of motion). 
By a suitable symplectic 
 close-to-the-identity
 nonlinear change of coordinates,
 the BNF
 puts the Hamiltonian in a simpler
 form in a  $\varepsilon$-neighborhood
 of the origin,
 for values of $\varepsilon$ sufficiently small. More precisely, after introducing
 action-angle variables,
 in the nonresonant case,
 the truncated system at order four
 is \textsl{integrated}, meaning
 its Hamiltonian depends 
 only on the actions,
 which are constants of motion,
 and not on the angles.
  The nonlinear frequencies
  are, then, easily evaluated
  as the derivatives of the Hamiltonian,
  i.e. the energy, with respect
  to the two actions.
  The resonant case is more
  complicated because
  the Hamiltonian also depends
   on one ``slow'' angle and 
 will be studied in detail elsewhere \cite{DL2}.
 The  key difference between resonant 
 and nonresonant regimes
relates to the mutual smallness of 
 the perturbative parameter $\varepsilon$
 and the detuning parameter
 $\sigma:=
 \omega_+-3 \omega_-$.
Specifically, we evaluate a constant
 $c_1$ such that we construct: 
\begin{enumerate}\label{zoomarine}
\item[-] a nonresonant  BNF, if $\varepsilon\leq c_1 \sqrt{|\sigma|}$
\label{mutua}
\item[-]  a resonant  BNF, if $ c_1 \sqrt{|\sigma|}<\varepsilon$
\end{enumerate}
(and $\varepsilon$ is small enough).
   A significant advantage of  Hamiltonian perturbation techniques is that
  one can estimate the magnitude of the maximal actions, equivalently the amplitudes, namely $\varepsilon$
  for which the procedure converges and, in addition, the magnitude of the 
  remainder in the Hamiltonian development, which contains monomials of order 6 or higher in  Cartesian variables (5 in the equations of motion).

\subsection*{Summary of the paper}

\subsubsection*{Application to the honeycomb metamaterial}
\vspace{-0.2cm}

In Section \ref{sec:2}
we discuss how the analysis presented in the following sections applies to the honeycomb metamaterial described earlier. Specifically, we highlight the set of resonant parameters and focus on evaluating the nonlinear bandgap in the nonresonant case.

\subsubsection*{Resonant and nonresonant Birkhoff Normal Form, and
applicability threshold}
\vspace{-0.2cm}

In Section \ref{sec:HamStr}
we reinterpret  the problem as a Hamiltonian 
system (see \eqref{ham}).
In Section \ref{Sec:BNF}
 we put the system, close to the origin, in  Birkhoff Normal Form both in the resonant and in the nonresonant case emphasizing the differences between them. 
 We then examine  the Hamiltonian truncated
 to within the fourth order, 
which is equivalent to the third order in the 
 equations of motion 
 (see \eqref{mozart} and \eqref{H4RES}
 for the nonresonant, and resonant cases, respectively), 
 as it captures the essential
 characteristics 
  of the overall motion.
By introducing action-angle variables,
 it becomes evident that 
  in the nonresonant scenario,  the truncated system is integrable, because the
 Hamiltonian 
depends solely 
  on the actions
 (see the first line in \eqref{secular}).
In contrast, in the resonant case,
 the truncated Hamiltonian, after a suitable linear change 
 of variables,
also  depends on one angle,  referred to as the  ``slow''
 angle 
(see the second line of \eqref{secular}), 
as its associated  
frequency is small or even zero at the exact resonance condition.
We compute the applicability threshold 
for the amplitudes in both the resonant and nonresonant cases of the Birkhoff Normal Form, see Figure \ref{Soglia_aPiuParRIS}.

\subsubsection*{The nonresonant case: 
 asymptotic solutions, estimate on the remainder and nonlinear frequencies
}
\vspace{-0.2cm}

Since the truncated
 Hamiltonian depends only on the actions,
  its nonlinear frequencies 
 are determined by the derivatives of the energy
 with respect to the actions (see 
 \eqref{omega_nonlineare})
 and all its motions are quasi-periodic
 (see Remark \ref{rem:KAM}).
The orbits of the complete system 
are the sums of  these quasi-periodic motions
(see 
\eqref{splendente})
plus a small remainder.
While the asymptotics of the solutions and the amplitude-dependent dispersion functions have already been computed using other methods, such as the multiple scales approach (see \cite{SW23mssp}), we recalculated them using Hamiltonian techniques and, additionally, provided explicit estimates on the remainders over an extended time interval (see \eqref{atene}).

\section{Application to the honeycomb metamaterial bandgap}\label{sec:2}

\subsection{Linear bandgap}

The (linear) bandgap is defined as
 the interval
\begin{equation}\label{puglia}
\left[\max_{(\tilde{k}_1,\tilde{k}_2)\in \triangle}\omega_-(\tilde M, \tilde K, \tilde{k}_1,\tilde{k}_2),\min_{(\tilde{k}_1,\tilde{k}_2)\in\triangle}\omega_+(\tilde M, \tilde K, \tilde{k}_1,\tilde{k}_2)\right]\,,
\end{equation}
	whose width is
\begin{equation}\label{Wlineare}
W=W(\tilde M, \tilde K):=
\min_{(\tilde{k}_1,\tilde{k}_2)\in\triangle}\omega_+(\tilde M, \tilde K, \tilde{k}_1,\tilde{k}_2)-
\max_{(\tilde{k}_1,\tilde{k}_2)\in \triangle}\omega_-(\tilde M, \tilde K, \tilde{k}_1,\tilde{k}_2)
\,,
\end{equation}	
and $\triangle$ indicates the IBZ.		
Since the gradients of $\omega_-$
and $\omega_+$ (with respect to $(\tilde{k}_1,\tilde{k}_2)$)
	never vanish
	 in the interior of $\triangle$,
	maxima and minima are attained on 
the boundary	$\partial \triangle$;
in particular,
\begin{eqnarray}\label{mela}
&&\max_{(\tilde{k}_1,\tilde{k}_2)\in \triangle}\omega_-(\tilde M, \tilde K, \tilde{k}_1,\tilde{k}_2)=
\max_{(\tilde{k}_1,\tilde{k}_2)\in \partial\triangle}\omega_-(\tilde M, \tilde K, \tilde{k}_1,\tilde{k}_2)
=\omega_-(\tilde M, \tilde K, \textstyle\frac43\pi,0)\,,
\\
&&\min_{(\tilde{k}_1,\tilde{k}_2)\in\triangle}\omega_+(\tilde M, \tilde K, \tilde{k}_1,\tilde{k}_2)=
\min_{(\tilde{k}_1,\tilde{k}_2)\in\partial\triangle}
\omega_+(\tilde M, \tilde K, \tilde{k}_1,\tilde{k}_2)
=\omega_+(\tilde M, \tilde K, 0,0)\,.
\nonumber
\end{eqnarray}
Hence, for every pair $(\tilde M,\tilde K)$, the maximum of the linear acoustic frequency 
is attained at  ${\bf X}$, while the minimum of
the linear optical frequency 
is attained at  ${\bf \Gamma}$.
We anticipate that, in evaluating the {\sl nonlinear} bandgap,  
${\bf X}$ will play a crucial role, more than ${\bf \Gamma}$.
As we will see,  the shift  of the maximum
of the acoustic frequency due to the nonlinearity is more significant than that 
of the minimum of the optical frequency.

\subsection{Resonant parameters}

In light of the  importance of the point ${\bf X}$
for the nonlinear bandgap,
we define an {\sl ${\bf X}$-resonant} pair as a pair $(\tilde M,\tilde K)$ such that 
the linear acoustic and optical frequencies evaluated at $(\tilde k_1,\tilde k_2)={\bf X}=(\frac43\pi,0)$
are in a 3:1 resonance, namely,
\begin{equation}\label{casetta}
\omega_{+}(\tilde M, \tilde K, \textstyle\frac43\pi,0)=
3\omega_{-}(\tilde M, \tilde K, \textstyle\frac43\pi,0)\,.
\end{equation}
The set of ${\bf X}$-resonant pairs is formed by two (red) curves plotted in Figure \ref{buoni}.
Since, here and henceforth, we will consider only  pairs in the
reference rectangle 
$[0.05, 0.3]\times[1, 20]$ (see Figure \ref{buoni}),
the set of ${\bf X}$-resonant pairs reduces to a unique (piece of) curve, which we refer to as
$\mathfrak R$,
specifically the solid red curve in Figure \ref{buoni}.

\begin{figure}[h!]
		\center	
	\includegraphics[width=16cm,height=12cm,keepaspectratio]{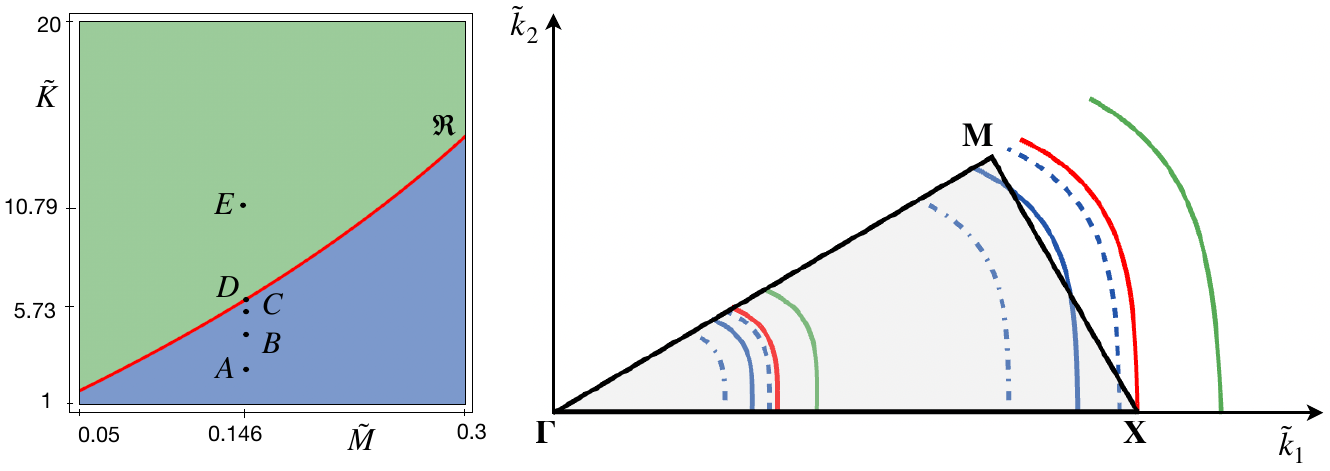} 
			\caption{(Left) The reference rectangle
			$[0.05,0.3]\times[1,20]$. 
			For every pair $(\tilde M,\tilde K)$
			in the blue region, there exist two resonant curves intersecting $\triangle$, while for pairs in the green region, there exists only one resonant curve intersecting $\triangle$.
						The red curve that separates the two regions
			is exactly the ${\bf X}$-resonant curve
			$\mathfrak R$. \\
(Right) Resonant curves
in the $(\tilde{k}_1,\tilde{k}_2)$-plane	 for different values of the pairs $(\tilde M,\tilde K)$ and their intersection with
the boundary of the Brillouin triangle $\triangle$.
The six blue curves  correspond to three
different
points $(\tilde M,\tilde K)$ in the blue region of
the rectangle.
 Specifically, when 
$(\tilde M,\tilde K)=(0.146,2)=A$
the corresponding two blue curves, the
dotted-dashed ones, 
 have 4 intersections. 
When 
$(\tilde M,\tilde K)=(0.146,3.6)=B$
the corresponding two blue curves, the
solid ones, 
 have 6 intersections.
 When 
$(\tilde M,\tilde K)=(0.146,5)=C$
the corresponding two blue curves, the
dashed ones, 
 have 4 intersections.
 The red curves, corresponding to
$(\tilde M,\tilde K)=(0.146,5.73)=D\in\mathfrak R$, have 3 intersections.  Finally the green
curves, corresponding to
$(\tilde M,\tilde K)=(0.146,10.79)=E$, have only 2 intersections. 
			}			\label{sogno}
		\end{figure}			 

Given a pair $(\tilde M,\tilde K)$ 
we call a set in the $(\tilde{k}_1,\tilde{k}_2)$-plane as {\sl resonant}
if every point in the set satisfies
the 3:1 resonance condition
 $3\omega_-
 (\tilde M,\tilde K,\tilde{k}_1,\tilde{k}_2)
 =\omega_+ (\tilde M,\tilde K,\tilde{k}_1,\tilde{k}_2)$.
For
any pair
$(\tilde M,\tilde K)$ within the rectangle 
$[0.05, 0.3]\times[1, 20]$ (see Figure \ref{sogno} (left))
 there are always one or two 
\textsl{resonant curves} in the
 $(\tilde{k}_1,\tilde{k}_2)$-plane,
 that  intersect
 the Brillouin triangle
$\triangle$ (see Figure \ref{sogno} (right)).
The  ${\bf X}$-resonant curve $\mathfrak R$
divides the rectangle $[0.05, 0.3]\times[1, 20]$ into two distinct regions:  the region above and the region below $\mathfrak R$, corresponding to
 the green region and the blue region in Figure
\ref{sogno} (left), respectively. 
For any fixed pair $(\tilde M,\tilde K)$ in the green region, there is only one 
resonant curve
in the plane of wave numbers $(\tilde k_1,\tilde k_2)$, that intersects the Brillouin triangle
(the green curve in Figure \ref{sogno} (right)).
Conversely, for every fixed pair $(\tilde M,\tilde K)$ in the blue region, there  are two 
resonant curves
in the plane of wave numbers $(\tilde k_1,\tilde k_2)$, that intersect the Brillouin triangle
(shown as the blue curves in Figure \ref{sogno} (right)).
Finally, in the limit case when the pair $(\tilde M,\tilde K)$ lies on the
${\bf X}$-resonant curve $\mathfrak R$, two 
resonant curves intersect the Brillouin triangle
in the $(\tilde k_1,\tilde k_2)$-plane, but one of them intersects
$\triangle$ only at the point ${\bf X}$
(see the red curves in Figure \ref{sogno} (right)).

\subsection{Admissible amplitudes}\label{sec:admissible}

\begin{figure}[h!]
\center
\includegraphics[width=14cm,height=14cm,keepaspectratio]{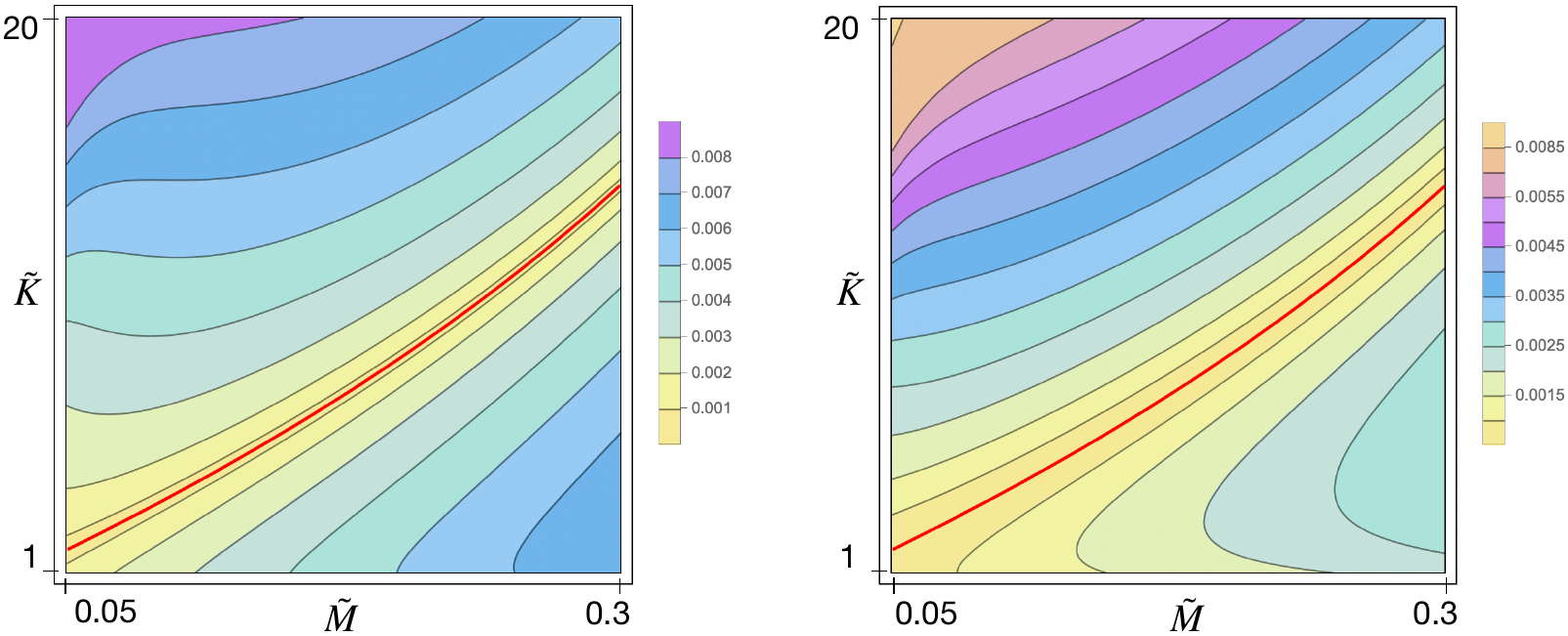}
\caption{Contour plots of the admissible 
wave amplitudes of the
optical (left) and acoustic (right) modes
in the plane $(\tilde M,\tilde K)$.
Here $N_3=-10^4$
and the wave numbers are chosen
in $\mathbf X$. Note that the value
decreases to zero as it approaches  the ${\bf X}$-resonant curve
 $\mathfrak R$.}
\label{Soglia_aPiuPar}
\end{figure}

\begin{figure}[h!]
			\centering			
			\includegraphics[width=7.5cm,height=7cm,keepaspectratio]{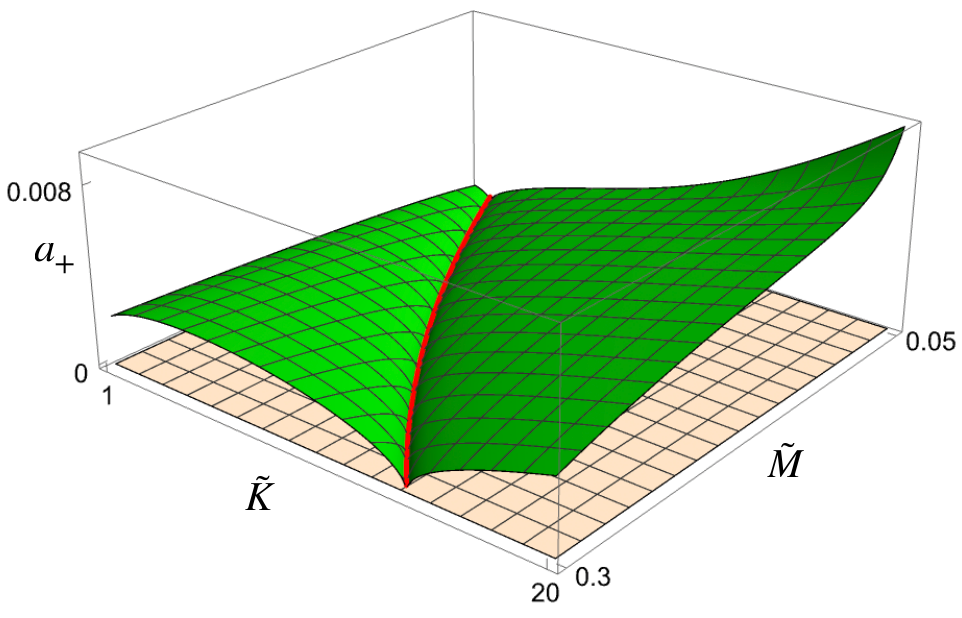} \ \ \ \
			\includegraphics[width=7.5cm,height=7cm,keepaspectratio]{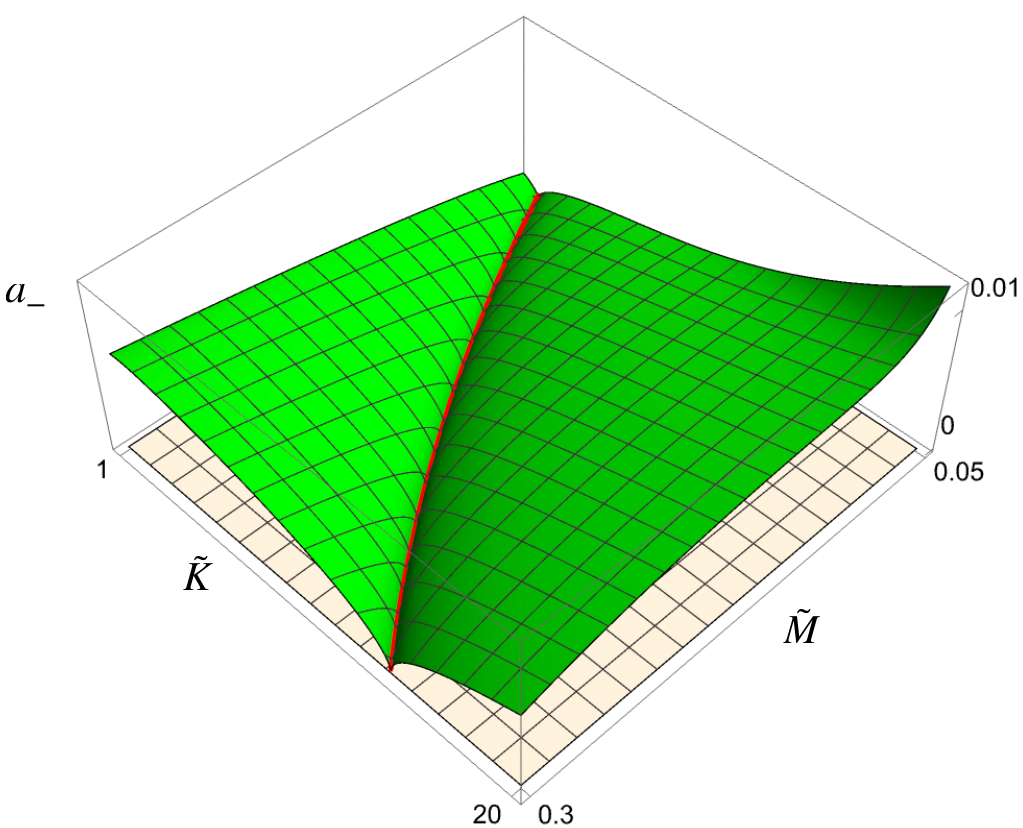}
			\caption{
		Admissible 
wave amplitudes of the
optical (left) and acoustic (right) modes are plotted
as functions of $\tilde M,\tilde K$.
The $(\tilde M,\tilde K)$-plane is plotted for 
comparison.
The ${\bf X}$-resonant curve  $\mathfrak R$
is shown in red.
			Here  $N_3=- 10^4$ (softening).
			}\label{a_Ticks}
		\end{figure}

\begin{figure}[h!]
\center
\includegraphics[width=12cm,height=12cm,keepaspectratio]{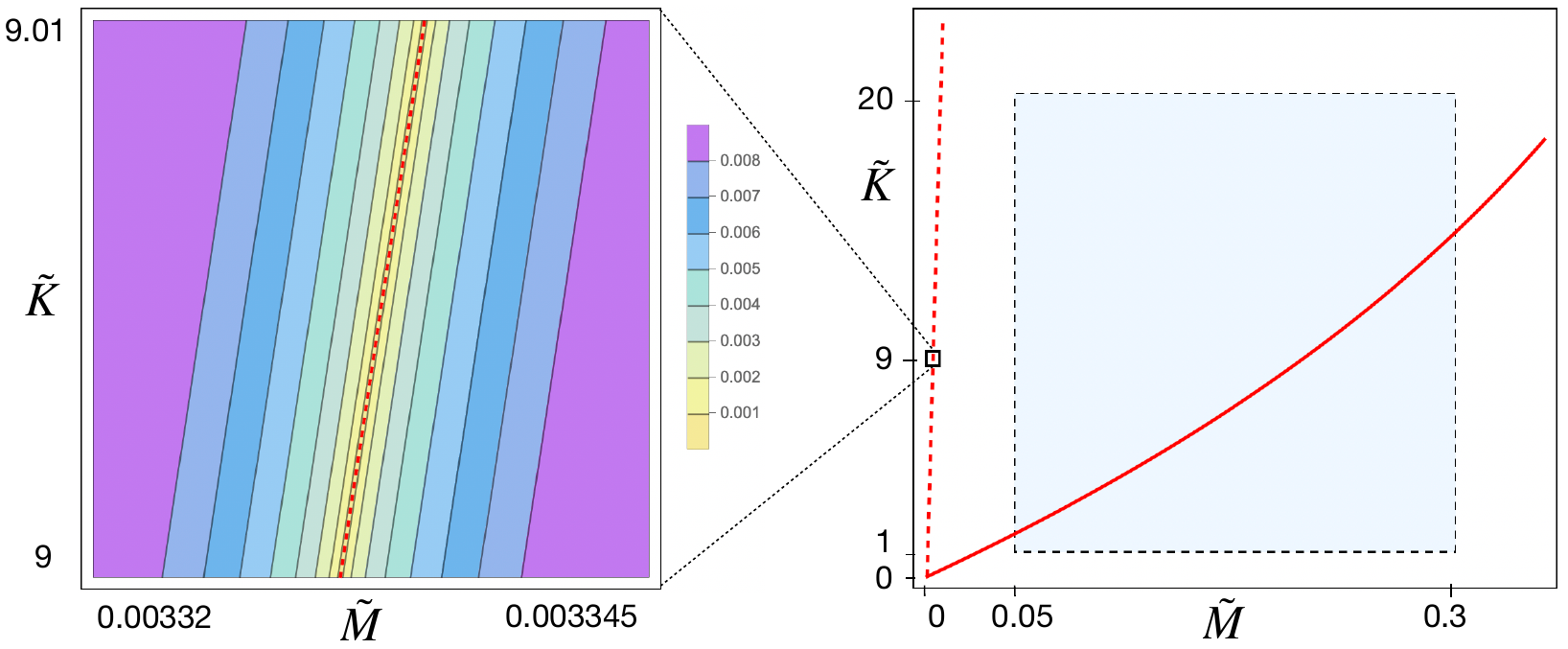}
\caption{Contour plot of the admissible 
wave amplitudes of the
 acoustic  mode (left)
in a small rectangle in the $(\tilde M,\tilde K)$-plane 
containing a small part of the ${\bf X}$-resonant dashed red curve which lies outside the light blue reference rectangle 
$[0.05,0.3]\times[1,20]$
(right).
Here $N_3=-10^4$
and the wave numbers are chosen
at $\mathbf X$. Note that the amplitude
decreases to zero as it approaches the ${\bf X}$-resonant curve.}
\label{rossa2}
\end{figure}

\begin{figure}[h!]
\center
\includegraphics[width=16cm,height=16cm,keepaspectratio]{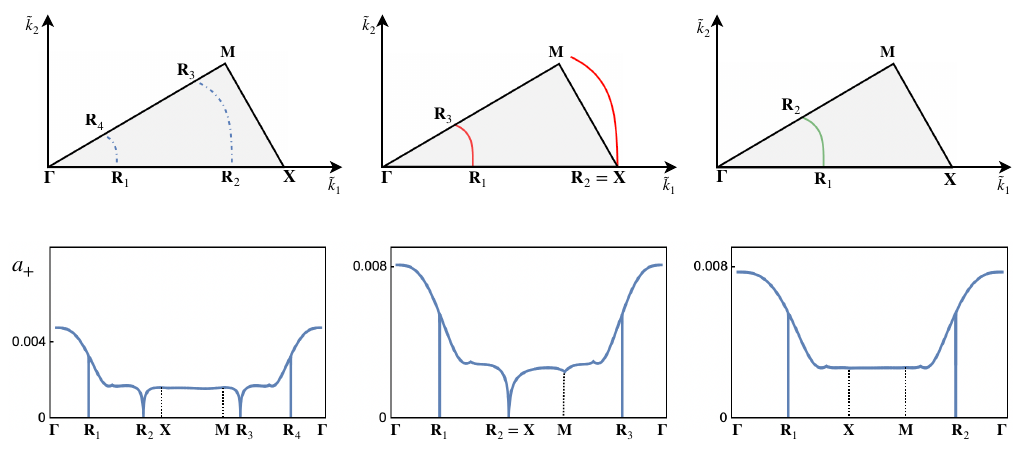}
\caption{
(Top)
Intersection points, 
${\bf R}_i$, between the
resonant curves
in the $(\tilde{k}_1,\tilde{k}_2)$-plane and
the boundary of the Brillouin triangle 
 for $(\tilde M, \tilde K)=(0.146,2)$
 (left), $(\tilde M, \tilde K)=(0.146,5.73)$
 (center),
 $(\tilde M, \tilde K)=(0.09,8)$
 (right).
 Compare with Figure 
\ref{sogno} (right).
\\
(Bottom)
Admissible 
wave amplitude, $a_+$,  of the optical mode
as a function of the wave numbers 
 on the boundary of the Brillouin triangle 
 for the same values of $(\tilde M, \tilde K)$ as above.
    Here, $N_3=-10^4$.
 Note that the value
decreases to zero on
the corresponding  four (left),    three (center),
two (right), resonant points  ${\bf R}_i$.
}\label{aPiuPar2}
\end{figure}

If we assume that we are away from the 3:1
resonance and that 
$N_3(a_-^2+a_+^2)$ is sufficiently small,
$a_\pm$ being the amplitudes of the optical and acoustic wave modes, respectively,
and $N_3$ is the strength of the nonlinearity, the formulas for the 
nonlinear dispersion curves 
$\omega_\pm^{\rm nl}=
\omega_\pm^{\rm nl}
(\tilde M, \tilde K, \tilde{k}_1,\tilde{k}_2,N_3,a_-, a_+)$
are\footnote{$\phi_i^\pm$, $i=1,2$, are the entries of the modal matrix, see \eqref{simultaneo}. They depend on 
$\tilde M$, $\tilde K$, $\tilde k_1$ and $\tilde k_2$.}
\begin{eqnarray}\label{formulaSL1}
			\omega_-^{\rm nl}
&=&		
\omega_-+N_3\left(
	\frac{3}{8\omega_-}
				(\phi_2^-)^4 a_-^2
				+
				\frac3{4\omega_-}
			(\phi_2^-)^2(\phi_2^+)^2
				a_+^2
				\right)
				\,,
				\nonumber
					\\
			\omega_+^{\rm nl}
			&=&
			\omega_+ +N_3\left(
	\frac{3}{8\omega_+}
				(\phi_2^+)^4
				a_+^2
				+
				\frac3{4 \omega_+}(\phi_2^-)^2(\phi_2^+)^2 a_-^2
				\right)
	\end{eqnarray}
(see formula (32) of \cite{SW23jsv}).
Here we discuss the validity of formula \eqref{formulaSL1}
in relation to the influence of the internal resonances,
once $N_3$ is fixed, e.g. $N_3=- 10^4$ (softening) or $10^4$ (hardening).
We first note that, since formula \eqref{formulaSL1} is a truncation 
at the second order of the expansions
of $\omega_\pm^{\rm nl}$ in series of $a_\pm$,
its validity is ensured only for sufficiently small $a_\pm$.
In Subsection \ref{subsec:BNF} we will
provide a (general) estimate (see \eqref{ampiezzeiniziali})
on the magnitude of $a_+$ and $a_-$
for which the above expansions are valid.
An important point to note is that this bound
strongly depends on the value of the resonance detuning expressed by
$|3 \omega_--\omega_+|$.
Specifically, the smaller $|3 \omega_--\omega_+|$, the smaller the allowed values of $a_+$ and $a_-$ must be.
When the linear frequencies are in a 3:1 
ratio, that is, $3 \omega_--\omega_+=0$,
$a_+$ and $a_-$ vanish. 
In particular, for $(\tilde k_1,\tilde k_2)=\mathbf X$,
the admissible amplitudes vanish for parameter values $(\tilde M, \tilde K)$
within the ${\bf X}$-resonant curve $\mathfrak R$.
The (maximum) admissible values of $a_+$ and $a_-$
as functions of $(\tilde M, \tilde K)$
with fixed $(\tilde k_1,\tilde k_2)=\mathbf X$
can be found in Figures \ref{Soglia_aPiuPar},
 \ref{a_Ticks}, and \eqref{aPiuPar2}.
 Incidentally, we note that a similar phenomenon 
 appears also around the other ${\bf X}$-resonant curve,
 namely 
the dashed red one in Figure 
\eqref{rossa2}, which, however, is outside 
the reference rectangle being considered.

By changing perspective,  we can fix $(\tilde M, \tilde K)$
and examine the variation of $a_\pm$  along  $\partial\triangle$.
 Note that the value of $a_\pm$
decreases to zero at certain resonant points. 
      The number of such resonant points varies  between 2 and 6
 depending on $(\tilde M, \tilde K)$. Let us refer to Figure \ref{sogno} (right).
 If the pair $(\tilde M, \tilde K)$ belongs to the green region
 of Figure \ref{sogno} (left), there is only one green curve intersecting  $\triangle$ and crossing   $\partial\triangle$ at two points.
When the pair $(\tilde M, \tilde K)$ belongs to $\mathfrak R$, namely the solid red curve 
 of Figure \ref{sogno} (left),  there are two curves intersecting  $\triangle$:
 the one on the right crosses  $\partial\triangle$ only at $\mathbf X$, while the one on the 
 left crosses  $\partial\triangle$ at two points.
 Finally, if the pair $(\tilde M, \tilde K)$ belongs to the blue region
 of Figure \ref{sogno} (left)  there are two curves intersecting  $\triangle$:
 the one on the 
 left crosses  $\partial\triangle$ at two points but 
 the one on the right may cross  $\partial\triangle$ only
 at two, three\footnote{This case, that occurs when the curve
 crosses $\partial\triangle$ at $\mathbf M$, is not shown in 
 Figure \ref{sogno} (right).} or four points.

The values of the admissible wave
amplitude $a_+$ 
on $\partial\triangle$ are shown 
 in Figure \ref{aPiuPar2} for three different pairs $(\tilde M, \tilde K)$.

In conclusion, we can say that, due to the presence of the
3:1 resonance, formula \eqref{formulaSL1} becomes invalid
in a neighborood   of the points ${\bf R}_i$,
when the parameters are resonant, namely where they
give rise to an exact 3:1 resonance between 
acoustic and optical  frequencies. 
In this case the expression of the nonlinear frequencies  
$\omega_\pm^{\rm nl}$ is much more complicated, due to the internal resonance,
and will be given elsewhere (see \cite{DL2}).
For a discussion on these issues see Remark
\ref{remo} below.

\subsection{Nonlinear bandgap}\label{sec:provola}

The nonlinear bandgap is the interval 
given by 
\begin{equation}\label{grecia}
\left[\max_{(\tilde{k}_1,\tilde{k}_2)\in \partial\triangle}
\omega_-^{\rm nl},
\ \ \min_{(\tilde{k}_1,\tilde{k}_2)\in\partial\triangle}
\omega_+^{\rm nl}
\right]\,,
\end{equation}
whose width is
\begin{equation}\label{molise}
W^{\rm nl}=W^{\rm nl}(\tilde M, \tilde K):=
\min_{(\tilde{k}_1,\tilde{k}_2)\in\partial\triangle}\omega_+^{\rm nl}(\tilde M, \tilde K, \tilde{k}_1,\tilde{k}_2)-
\max_{(\tilde{k}_1,\tilde{k}_2)\in \partial\triangle}\omega_-^{\rm nl}(\tilde M, \tilde K, \tilde{k}_1,\tilde{k}_2)
\,.
\end{equation}	
Note that the nonlinear bandgap and its width are, as in the linear case,
  functions of the pair $(\tilde M,\tilde K)$.
		
\begin{figure}[h!]
			\centering			
			\includegraphics[width=13cm,height=13
		cm,keepaspectratio]{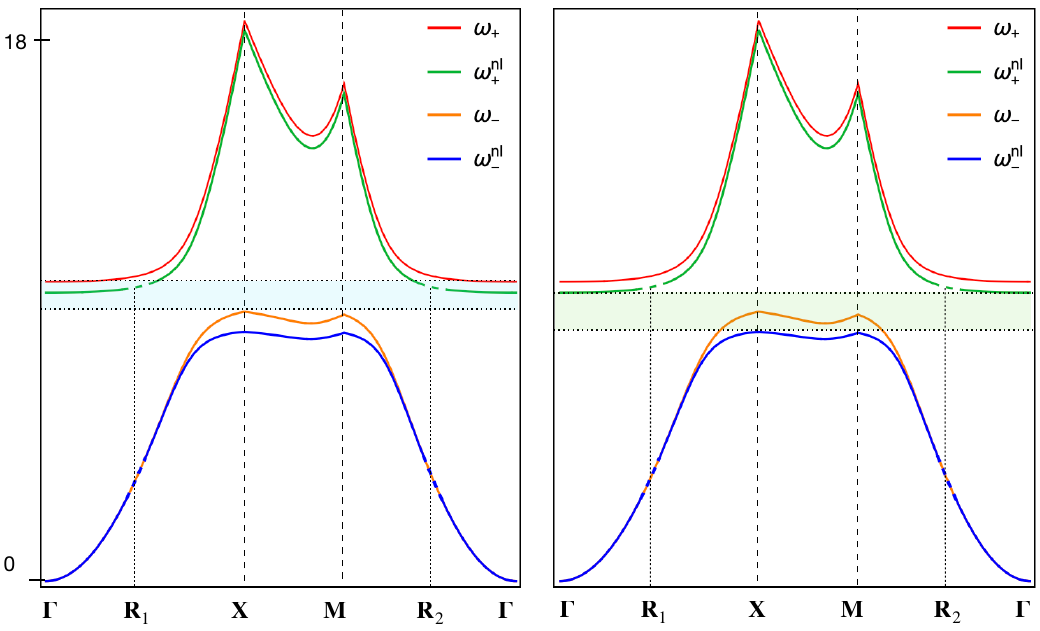}
			\caption{(Softening case) Comparison between the linear (on the left) and nonlinear (on the right) bandgap computed 
			for 
			 $\tilde M=0.09$, $\tilde K=8$, $N_3=-10^4$. According to Figure
			 \ref{Soglia_aPiuPar}
			 the corresponding admissible
			 initial amplitudes are $a_-=0.0036$
			and $a_+=0.0025$.
			$\omega_+$ and $\omega_-$
			are plotted in red and orange, respectively.
			The nonlinear frequencies 
			 $\omega_+^{\rm nl}$ and 
			$\omega_-^{\rm nl}$ given by formula \ref{formulaSL1}
			are plotted in green and blue, respectively. 
			Note that these 
 two curves are dashed, corresponding to the two neighborhoods of the resonant points ${\bf R}_1$ and ${\bf R}_2$. In these intervals,
 formula \ref{formulaSL1} does not hold and the effective  
 $\omega_+^{\rm nl}$ and 
			$\omega_-^{\rm nl}$ are slightly different (they are not shown here). Nonetheless, the minimum of $\omega_+^{\rm nl}$ and the maximum of
			$\omega_-^{\rm nl}$
			are still attained at ${\bf \Gamma}$ and ${\bf X}$, respectively,
			as for $\omega_+$ and $\omega_-$.
			Note that, for such values of the pair $(\tilde M,\tilde K)$,
			a softening nonlinearity gives rise to an increment of the bandgap width.
			 }
			 \label{piscina1}
		\end{figure}

\begin{figure}[h!]
			\centering			
			\includegraphics[width=13cm,height=13cm,keepaspectratio]{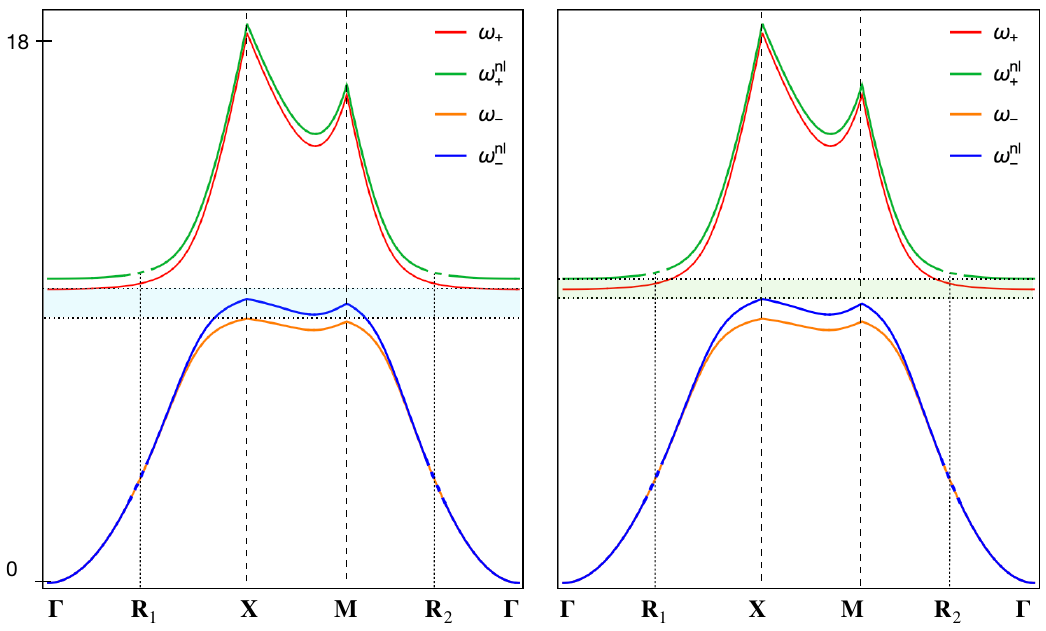}
			\caption{(Hardening case) Comparison between the linear (left) and nonlinear (right) bandgap computed 
			for 
			 $\tilde M=0.09$, $\tilde K=8$, $N_3=10^4$. According to Figure
			 \ref{Soglia_aPiuPar}
			 the corresponding admissible
			 wave amplitudes are $a_-=0.0036$
			and $a_+=0.0025$.
			The same considerations of Figure \ref{piscina1} hold.
			Note that for such values of the parameters
			a hardening nonlinearity gives rise to a shrinking of the bandgap.
			 }
			 \label{piscina2}
		\end{figure}

Given a pair $(\tilde M,\tilde K)$, the calculation of the nonlinear bandgap
 in \eqref{grecia}  
 presents no difficulties
 when the maximum and the minimum are attained away from the resonant
 points, where  formula \eqref{formulaSL1} holds.
 An example of this situation is shown in Figures 
 \ref{piscina1} and \ref{piscina2} in the softening ($N_3<0$) and hardening
 ($N_3>0$) case,
 respectively.
 

\begin{figure}[h!]
			\centering			
			\includegraphics[width=13cm,height=13
		cm,keepaspectratio]{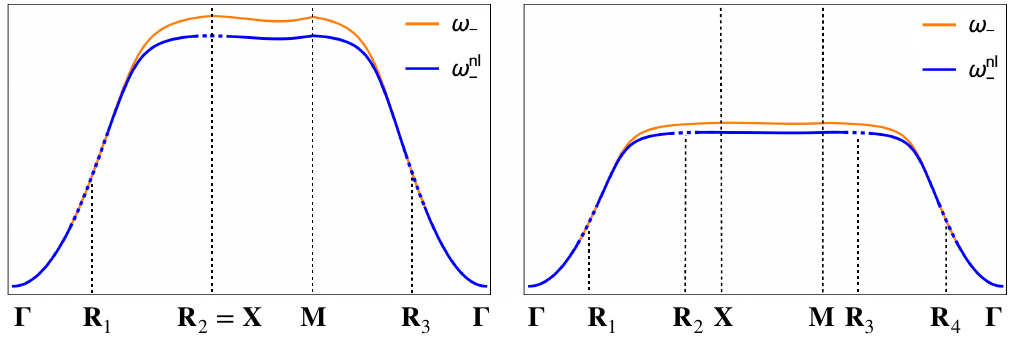}
			\caption{Linear and nonlinear acoustic frequencies
			for $N_3=-10^4$ and $(\tilde M,\tilde K)=(0.146, 5.73)$ 
			 (left)
			and $(\tilde M,\tilde K)=(0.146, 2)$ (right).
			 According to Figure
			 \ref{Soglia_aPiuPar}
			 the corresponding admissible
			 initial amplitudes are $a_-=0.0037$, $a_+=0.0022$
			 and $a_-=0.0019$, $a_+=0.0011$, respectively.
			Note that the nonlinear acoustic frequencies are dashed, corresponding to the  neighborhoods of the resonant points ${\bf R}_i$  since, in such intervals,
 formula \ref{formulaSL1} does not hold and the effective  
			$\omega_-^{\rm nl}$ is slightly different. 
			They are not shown here (see \cite{DL2} for details) but are closer to the linear 
			acoustic frequencies. This implies that the maximum
			of $\omega_-^{\rm nl}$ is attained
			in a neighborood of ${\bf R}_2$
			(both left and right).
			 }
			 \label{casibrutti}
		\end{figure}

 Note that, since we are considering pairs $(\tilde M,\tilde K)$ 
 belonging to  the rectangle $[0.05, 0.3]\times[1, 20]$,
 the IBZ vertex ${\bf \Gamma}$, where the minimum of the linear optical
 frequency is attained, is always far from being resonant.  
Therefore, henceforth, we will focus on 
the maximum of the acoustic frequency  because
it undergoes the most significant shifts
and may be resonant.
It turns out that, for the calculation of the nonlinear bandgap,
there are two  issues: i) when  $\mathbf X$ is resonant or nearly resonant,
ii) when $\omega_-$ has an almost flat maximum (at ${\bf X}$), so that, even if
${\bf X}$ is away from resonance,  $\omega_-^{\rm nl}$,
despite being close to  $\omega_-$, may attain its maximum at 
some resonant (or nearly resonant) point away from ${\bf X}$.
These two cases are shown in Figure \ref{casibrutti}.

The set of ``good'' pairs $(\tilde M,\tilde K)$ 
in  the rectangle $[0.05, 0.3]\times[1, 20]$, for which
 the maximum of $\omega_-^{\rm nl}$  and the minimum
 of $\omega_+^{\rm nl}$ are attained away from resonant
 points,
  is shown in light yellow in Figure~\ref{buoni}. 
More precisely the minimum of $\omega_+^{\rm nl}$  and the maximum
 of $\omega_-^{\rm nl}$ are attained at
 ${\bf \Gamma}$ and ${\bf X}$, respectively, with 
 ${\bf X}$ away from being resonant.
 In this case the  nonlinear bandgap width  can be evaluated  by  formula \eqref{formulaSL1} as
\begin{equation}\label{abruzzobis}
W^{\rm nl}(\tilde M, \tilde K)=
\omega_+^{\rm nl}(\tilde M, \tilde K, 0,0)-
\omega_-^{\rm nl}(\tilde M, \tilde K, \textstyle\frac43\pi,0)
\,.
\end{equation}

Note that in order to maximize the nonlinear bandgap width, as given by formula
\eqref{formulaSL1}, it is convenient to  take 
the wave amplitudes as large as possible.
Therefore, in light of the above discussion,
we choose the admissible amplitudes described in 
Subsection \ref{sec:admissible}. 

Recalling \eqref{Wlineare}, let us introduce the 
maximum percentage difference in the bandgap
as
\begin{equation}\label{Bper}
B_{\rm per}=B_{\rm per}(\tilde M,\tilde K):=100\left(
\frac{W^{\rm nl}(\tilde M, \tilde K)}{W(\tilde M, \tilde K)}-1
\right)\,.
\end{equation}
The level curves of $B_{\rm per}$ 
are plotted in Figure \ref{ritz}.
In creating this plot, for the light yellow zone shown in Figure
\ref{buoni}, we can use the representation formula \eqref{formulaSL1}
in addition to \eqref{Wlineare}, \eqref{abruzzobis} and \eqref{Bper}.
On the other hand, in the complementary light purple zones in 
Figure
\ref{buoni} we cannot use \eqref{formulaSL1} since, as we have already
noted,  for values of the pair $(\tilde M,\tilde K)$ in the light purple region,
 the maximum is attained near a resonant point on $\partial\triangle$.
Then, to evaluate the nonlinear bandgap in the light purple
 zone, one must find a new ``resonant'' representation
 formula for the nonlinear frequencies.
 Such formula is quite complex and  is derived through the resonant BNF, as  
 detailed  in \cite{DL2}, 
 where we also discuss   in details how to derive Figure \ref{ritz}.

\section{The Hamiltonian structure}\label{sec:HamStr}

In this section, after
introducing the optical and acoustic 
modes,
we identify the system
in \eqref{autostrada}  as  Hamiltonian,
see \eqref{ham} below, and we evaluate
the coefficients of the Hamiltonian, see
\eqref{fij}.
Set
$$
\mathbf \Lambda:=
 \left(
\begin{array}{cc}
	\omega_-^2&0\\
       0& \omega_+^2
	\end{array}	
\right)\,, 
$$
where
$\omega_-^2<\omega_+^2$ are the positive
eigenvalues of $\mathtt M^{-1}\mathtt K$ and
 $0<\omega_-<\omega_+$.
Since $\mathtt M$ is symmetric and 
$\mathtt K$ is diagonal,
there exists a $2\times 2$ matrix $\mathbf  \Phi$  such that
\begin{equation}\label{simultaneo}
\mathbf \Phi^T \mathtt M\mathbf \Phi=\mathbf  I\,,
\qquad
\mathbf \Phi^T \mathtt K \mathbf \Phi=\mathbf \Lambda\,,
\qquad
\mathbf  \Phi
 =\left(
\begin{array}{cc}
	\phi_1^-& \phi_1^+\\
       \phi_2^-& \phi_2^+
	\end{array}	
\right)\,,
\end{equation}
where $\mathbf  I$ is the identity matrix. 
Consider the change of variables
$$\left(
		\begin{array}{c}
			v \\
			y
		\end{array}
		\right)
		=\mathbf \Phi \mathbf q\,,\qquad
		\mathbf q:=\left(
		\begin{array}{c}
			q_1 \\
		q_2
		\end{array}
		\right)\,.
		$$
	By  \eqref{simultaneo}
	the system in \eqref{autostrada}
	is transformed into 
	\begin{equation}\label{fettuccine}
	\ddot {\mathbf q}+\mathbf \Lambda \mathbf q=\mathbf n(\mathbf q)\,,\qquad
	\mathbf n(\mathbf  q)=\left(
		\begin{array}{c}
			n_1 \\
		n_2
		\end{array}
		\right):=-\mathbf \Phi^T\left(
		\begin{array}{c}
	M_3 (\phi_1^- q_1+\phi_1^+ q_2)^3 \\
	N_3 (\phi_2^- q_1+\phi_2^+ q_2)^3 \\
	\end{array}	
	\right)\,.
\end{equation}
In particular
$$
			\begin{array}{ll}
					n_1=-\phi_1^- M_3(\phi_1^- q_1+\phi_1^+ q_2)^3-\phi_2^- N_3(\phi_2^- q_1+\phi_2^+ q_2)^3\\
					n_2=
					-\phi_1^+M_3(\phi_1^- q_1+\phi_1^+ q_2)^3-\phi_2^+N_3(\phi_2^- q_1+\phi_2^+ q_2)^3\,.
			\end{array}		
	$$
Introducing the momenta 
$\dot {\mathbf q}=\mathbf p=\left(
		\begin{array}{c}
			p_1 \\
		p_2
		\end{array}
		\right)$,	 the system in \eqref{fettuccine}
	is Hamiltonian with Hamiltonian
	\begin{equation}\label{ham}
		H(\mathbf p,\mathbf q)=\frac 12 (p_1^2+p_2^2) + \frac12 \omega_-^2 q_{1}^2+ \frac12 \omega_+^2 q_{2}^2+
		f(\mathbf q)\,,
	\end{equation}
where 
	\begin{equation}
				f(\mathbf q)
			:= \frac 14 M_3 (\phi_1^- q_1+\phi_1^+ q_2) ^4+\frac 14 N_3 (\phi_2^- q_1+\phi_2^+ q_2)^4\,.
				\label{potential}
	\end{equation}
Indeed it is straightforward to see that  Hamilton's
equations
$\dot {\mathbf p}=-\partial_{\mathbf q} H$, $\dot {\mathbf q}=\partial_{\mathbf p} H=\mathbf p$
are equivalent to the system 
in \eqref{fettuccine}.
Since $f(\mathbf q)$ is a homogeneous polynomial
of degree $4$ we write
	\begin{equation}\label{fij}
	f(\mathbf q)=\sum_{i+j=4}f_{i,j}q_1^i q_2^j\,,
	\qquad
	\mbox{with}
	\quad			f_{i,j}:=\frac{6}{i!j!}
				\Big(
				(\phi_1^-)^i(\phi_1^+)^jM_3+(\phi_2^-)^i(\phi_2^+)^jN_3
				\Big)\,.
	\end{equation}
Introducing coordinates 
$\mathbf Q=(Q_1,Q_2)$, $\mathbf P=(P_1,P_2)$ through 
\begin{equation}\label{islanda}
			p_1=\sqrt{\omega_-}P_1\qquad p_2=\sqrt{\omega_+}P_2\qquad q_1=\frac1{\sqrt{\omega_-}}Q_1\qquad q_2=\frac1{\sqrt{\omega_+}}Q_2\,,
		\end{equation}
it turns out   that the Hamiltonian in the new variables
reads:
	\begin{equation}\label{hamPQ}
		\mathtt H(\mathbf P,\mathbf Q):=\omega_- \frac{P_1^2+Q_1^2}{2}+\omega_+ \frac{P_2^2+Q_2^2}{2}+ f\left(\frac{Q_1}{\sqrt{\omega_-}},\frac{Q_2}{\sqrt{\omega_+}}\right)\,.
	\end{equation}
By introducing complex coordinates, 
$\mathrm i=\sqrt{-1}\in \mathbb C $, $\mathbf z=(z_1,z_2)\in \mathbb C^2$
		\begin{equation}\label{islanda2}
			z_j=\frac{Q_j+\mathrm i P_j}{\sqrt 2}\quad \bar z_j=\frac{Q_j-\mathrm i P_j}{\sqrt 2}\qquad j=1,2
	\end{equation}
		the Hamiltonian reads:
\begin{equation}\label{hamiltonian}
		\mathtt H(\mathbf z,\bar{\mathbf z})
			=
				\mathtt N(\mathbf z,\bar{\mathbf  z})
				+
				{\mathtt G}(\mathbf z,\bar{\mathbf  z})
\end{equation}
where
\begin{equation}\label{GG}
\mathtt N(\mathbf z,\bar{\mathbf  z}):=
				\omega_- z_1\bar z_1+
				\omega_+ z_2\bar z_2\,,
				\qquad
{\mathtt G}(\mathbf z,\bar{\mathbf  z})
			:=
			f\left(\frac{z_1+\bar z_1}{\sqrt{2\omega_-}},\frac{z_2+\bar z_2}{\sqrt{2\omega_+}}\right)\,.
\end{equation}
Note that, in complex coordinates,
 Hamilton's equations become
	\begin{equation}\label{HE}
		 \dot z_j=-\mathrm{i} \partial_{\bar z_j}\mathtt H\,, \ \ \dot{\bar{z}}_j= \mathrm{i} \partial_{z_j}\mathtt H\,.
	\end{equation}
In the following
we employ the multi-index notation
\begin{equation}\label{fuga}
		P(\mathbf z,\bar{\mathbf  z})=\sum_{(\bal,\bbe)\in\mathbb N^2\times\mathbb N^2} P_{\bal,\bbe}
		\mathbf  z^\bal
		\bar{\mathbf z}^\bbe
		\end{equation}
for suitable coefficients $P_{\bal,\bbe}\in\mathbb C$
with $\mathbf z^\bal=z_1^{\alpha_1}z_2^{\alpha_2}$ (analogously
for $\bar{\mathbf z}^\bbe$). 
In this notation, recalling
\eqref{fij} and \eqref{GG},
 we rewrite ${\mathtt G}$ 
 as\footnote{Where, for integer vectors
 $\bal,\bbe$ we set
 $|\bal+\bbe|:=\alpha_1+\alpha_2+\beta_1+\beta_2$.
}
	\begin{equation}\label{G}
			{\mathtt G}(\mathbf z,\bar{\mathbf  z})
						=
				\sum_{i+j=4}\frac{f_{i,j}}{4(\sqrt{\omega_-})^{i}
				(\sqrt{\omega_+})^{j}}
				(z_1+\bar z_1)^i
				(z_2+\bar z_2)^j
			=
				\sum_{|\bal+\bbe|=4 }
				{\mathtt G}_{\bal,\bbe} \mathbf z^\bal \bar{\mathbf z}^\bbe
	\end{equation}
where\footnote{
Being
	\begin{equation}\label{binomio}
			(z_1+\bar z_1)^i(z_2+\bar z_2)^j
			=
		\sum_{h=0}^i 
			\left(
				\begin{array}{c}
					i\\
					h			
				\end{array}
			\right)
			z_1^h \bar z_1^{i-h}
		\sum_{k=0}^j
			\left(
				\begin{array}{c}
					j\\
					k			
				\end{array}
			\right)
			z_2^k \bar z_2^{j-k}
			=
\sum_{\alpha_1 +\beta_1=i,\  \alpha_2 +\beta_2=j} 
c_{\bal, \bbe}\mathbf z^{\bal}\bar{\mathbf z}^{\bbe}
	\end{equation}	
and, taking
$\alpha_1 +\beta_1=i,$ $\alpha_2 +\beta_2=j$,	
	 \begin{equation*}\label{bino}
	 \left(
				\begin{array}{c}
					i\\
					h			
				\end{array}
			\right)
			\left(
				\begin{array}{c}
					j\\
					k			
				\end{array}
			\right)=
			\left(
				\begin{array}{c}
					\alpha_1+\beta_1\\
					\alpha_1			
				\end{array}
			\right)
			\left(
				\begin{array}{c}
					\alpha_2+\beta_2\\
					\alpha_2			
				\end{array}
			\right)
			=
			\frac{(\alpha_1+\beta_1)!}{\alpha_1! \beta_1!}	
			\frac{(\alpha_2+\beta_2)!}{\alpha_2! \beta_2!}=c_{\bal,\bbe}\,.		
	 \end{equation*}} 
	\begin{equation}\label{Gab}
				{\mathtt G}_{\bal,\bbe}:=
				\frac{f_{\alpha_1+\beta_1,\alpha_2+ \beta_2}}
				{4(\sqrt{\omega_-})^{\alpha_1+\beta_1}
				(\sqrt{\omega_+})^{\alpha_2+\beta_2}}
				\frac{(\alpha_1+\beta_1)!}{\alpha_1!\beta_1!}
				\frac{(\alpha_2+\beta_2)!}{\alpha_2!\beta_2!}\,.
	\end{equation}
Note that ${\mathtt G}_{\bal,\bbe}=
{\mathtt G}_{\bbe,\bal}\in\mathbb R$.

\section{Resonances and Birkhoff Normal Forms}\label{Sec:BNF}

In this section
we use symplectic changes
of variables in order to put the Hamiltonian
in BNF both in the nonresonant and
in the resonant case, see
Theorem \ref{rossoscuro}.  
These changes of coordinates are 
constructed as Hamiltonian flows
of suitable ``generating functions''.
	
The aim of the BNF is to construct a symplectic change
of variables that ``simplifies'' the Hamiltonian
$\mathtt H$ in \eqref{hamiltonian}.
First note that a Hamiltonian $H$ depending 
only on 	$|z_1|^2$ and $|z_2|^2$
can be expressed as 
$H=\sum_{\bal} H_{\bal,\bal}
|\mathbf  z|^{2\bal}$
and is integrable; in particular $|z_1|^2$ and $|z_2|^2$
 are constants of motion
as one immediately recognises since, for $j=1,2$,
we have\footnote{Recall \eqref{HE}.}
	$$
		\frac{d}{dt}|z_j|^2=\frac{d}{dt}(z_j\bar z_j)=
		\dot z_j\bar z_j+z_j\dot{\bar z}_j=-\mathrm i \bar z_j\partial_{\bar z_j}
		H+\mathrm i z_j\partial_{z_j}
		H =0\,.
	$$
In light of the above considerations
 we aim to find out whether  it is possible to determine,
in a sufficiently small neighbourhood of the origin, a close-to-the-identity symplectic transformation that integrates 	
	$\mathtt H$ up to terms of degree 6 in $(\mathbf z,\bar{\mathbf  z})$,
	which are small near the origin.	
	This amounts to transforming 
	 $\mathtt H$ into
	${\mathtt N}+\bar{\mathtt H}_4+O(\|\mathbf z\|^6)$, with
	\begin{equation}\label{mozart}
		\bar{\mathtt H}_4:=
		\sum_{ |\bal|=2} {\mathtt G}_{\bal,\bal} |\mathbf  z|^{2\bal}
		=
		{\mathtt G}_{(2,0),(2,0)} 
		|z_1|^4
			+
		{\mathtt G}_{(1,1),(1,1)}
		|z_1|^2 |z_2|^2
		+
		{\mathtt G}_{(0,2),(0,2)} 
		|z_2|^4\,,
	\end{equation}
where, by recalling \eqref{Gab},
	\begin{eqnarray}
				{\mathtt G}_{(2,0),(2,0)} 
			&=& 
				\frac{3 f_{4,0}}{2 \omega_-^2}
				=\frac{3}{8\omega_-^2}
				\Big(
				(\phi_1^-)^4M_3+(\phi_2^-)^4N_3
				\Big)\,,
				\nonumber
				\\
				{\mathtt G}_{(1,1),(1,1)} &=&
				\frac{f_{2,2}}{\omega_- \omega_+}
				=\frac3{2\omega_- \omega_+}\Big(
				(\phi_1^-)^2(\phi_1^+)^2M_3+(\phi_2^-)^2(\phi_2^+)^2 N_3
				\Big)\,,
				\nonumber
				\\
				{\mathtt G}_{(0,2),(0,2)} 
			&=&
				\frac{3 f_{0,4}}{2 \omega_+^2}
				=\frac{3}{8\omega_+^2}
				\Big(
				(\phi_1^+)^4M_3+(\phi_2^+)^4N_3
				\Big)\,.
	\label{calippo}			
	\end{eqnarray}

\subsection{Symplectic change of coordinates}

Let 
$\left\{\cdot,\cdot\right\}$ be the Poisson brackets
defined as
	\begin{equation}\label{poisson}
			\{F,G\}:=\mathrm i \sum_{j=1,2}\big(\partial_{z_j}F\partial_{\bar z_j}G-\partial_{\bar z_j}F\partial_{z_j}G\big)\,.
	\end{equation}
Incidentally note that if $F$ and $G$ are homogeneous
polynomials of degree $d_F$ and $d_G$, respectively,
then
	\begin{equation}\label{paderborn}
		\{F,G\} \ \ \mbox{is a homogeneous polynomial of degree}
		\ \ d_F+d_G-2\,.
	\end{equation}

A typical way to produce symplectic transformation is by the time-$1$ Hamiltonian flow $\Phi^1_S$ generated by an auxiliary Hamiltonian  $S$, the ``generating function''.
Let us take as generating function $S$  a
homogeneous polynomial of degree 4.
If $S$ is small enough\footnote{The quantitative aspects will be detailed
in Subsection \ref{subsec:BNF}.},
the Hamiltonian in
\eqref{hamiltonian} is
 transformed  into
 $\mathtt H\circ\Phi^1_S$, which
  can be written as
  	\begin{equation}\label{bach}
			\mathtt H\circ\Phi_S^1
			=
{\mathtt N}\circ\Phi_S^1 +\mathtt G\circ\Phi_S^1
=
 {\mathtt N}
 +\{S,{\mathtt N}\}+\mathtt G
 +R\,,
 \end{equation}
 where, by \eqref{paderborn},
			$\mathtt G+\{S,{\mathtt N}\}$
			is a homogeneous polynomial of degree 4
			and
			the 
remainder\footnote{\label{footnote6} Take the Taylor
expansion at $t=0$ of the functions
$t\to {\mathtt N}\circ\Phi_S^t$
and $t\to {\mathtt G}\circ\Phi_S^t$
evaluated at $t=1$, of order
2 and 1, respectively.
Note that, by the complex
form of  Hamilton's equations,
 the Hamiltonian
flow $\Phi_S^t$ satisfies
$\frac{d}{dt}\Phi_S^t=
\mathrm i \mathbf J \nabla S(\Phi_S^t)$,
where $\mathbf J:= \left(
\begin{array}{cc}
	0&-\mathbf I_{2\times 2}\\
       \mathbf I_{2\times 2}& 0
	\end{array}	
\right)$
 is the standard
symplectic matrix, where
$\mathbf I_{2\times 2}$
is the $2\times 2$ identity matrix.
Recalling \eqref{poisson}, 
$\frac{d}{dt}({\mathtt G}\circ
\Phi_S^t)=(\nabla \mathtt G
)\circ
\Phi_S^t\, \frac{d}{dt}\Phi_S^t
=
\mathrm i (\nabla \mathtt G\mathbf J \nabla S)\circ
\Phi_S^t
=\{S,\mathtt G\}\circ
\Phi_S^t$.
Analogously 
$\frac{d^2}{dt^2}({\mathtt N}\circ
\Phi_S^t)=\{S,\{S,\mathtt N\}\}\circ
\Phi_S^t$.
}
 $$
R:=\int_0^1(1-t)\{S,\{S,{\mathtt N}\}\}
 \circ\Phi_S^t\, dt
+\int_0^1\{S,{\mathtt G}\}
 \circ\Phi_S^t\, dt\,.			
$$
		We look for a generating function 
		$S=\sum_{|\bal+\bbe|=4}S_{\bal,\bbe}\mathbf z^\bal \bar{\mathbf z}^\bbe$, which simplifies
		the term of order 4, namely $\mathtt G+\{S,{\mathtt N}\}$ (recall \eqref{G}).
		In particular, we search for a function   $S$ satisfying the ``homological equation'':
	\begin{equation}\label{beige}
			\mathtt G+\{S,{\mathtt N}\}=\bar{\mathtt H}_4\,,
	\end{equation}
	with $\bar{\mathtt H}_4$ defined in \eqref{mozart}
		so that by \eqref{bach} and \eqref{beige} we have
\begin{equation}\label{bach2}
			\mathtt H\circ\Phi_S^1
			=
{\mathtt N} +\bar{\mathtt H}_4+R
\end{equation}
with
 \begin{equation}\label{resto}
R=\int_0^1(1-t)\{S,
\bar{\mathtt H}_4-\mathtt G\}
 \circ\Phi_S^t\, dt
+\int_0^1\{S,{\mathtt G}\}
 \circ\Phi_S^t\, dt
 =
 \int_0^1\big\{S,
\bar{\mathtt H}_4+t (\mathtt G-\bar{\mathtt H}_4)\big\}
 \circ\Phi_S^t\, dt\,.			
	\end{equation}
The formal solution of \eqref{beige} is straightforward:
since, by \eqref{GG} and \eqref{poisson},
$$
\{S,{\mathtt N}\}=\mathrm i
\sum_{|\bal+\bbe|=4}
\bom\cdot(\bal-\bbe)
S_{\bal,\bbe}\mathbf z^\bal \bar{\mathbf z}^\bbe\,,
\qquad
\bom:=(\omega_-,\omega_+)\,,
$$
then the coefficients $S_{\bal,\bbe}$ can be chosen as 
\begin{equation}\label{vivaldi}
S_{\bal,\bbe}:=\mathrm i\frac{\mathtt G_{\bal,\bbe}}{\bom\cdot(\bal-\bbe)}\,,\qquad
\mbox{if}\ \ \ \bal\neq \bbe
\end{equation}
and $0$ otherwise.
It should be clear that the expressions in \eqref{vivaldi}
make sense only when the denominators
(``small divisors'') do not vanish.
This depends on the nonresonant properties
of the frequency vector
$\bom$.
 
 \subsection{Resonances}
 
We have to bound from below the 
quantities 
$|\bom\cdot(\bal-\bbe)|$ as 
$|\bal+\bbe|=4$ and $\bal\neq\bbe$.

\begin{pro}\label{dartagnan}
 If $\omega_+\neq 3 \omega_-$ then
 for every 
$(\bal,\bbe)\in\mathbb N^2\times\mathbb N^2$,
$|\bal+\bbe|=4$, $\bal\neq\bbe$, we have
 \begin{equation}\label{stoccolma}
 |\bom\cdot(\bal-\bbe)|\geq
\gamma:=\min\{\omega_-, \omega_+-\omega_-, |3 \omega_--\omega_+|\}>0\,.
\end{equation}
On the other hand, in the resonant case
 $\omega_+= 3 \omega_-$, 
  for every 
$|\bal+\bbe|=4$, $\bal\neq\bbe$,
$(\bal,\bbe)\neq (3,0,0,1), (0,1,3,0)$
 we have
 \begin{equation}\label{stoccolma2}
 |\bom\cdot(\bal-\bbe)|\geq
\gamma_{{\rm res}}:=\min\{\omega_-, \omega_+-\omega_-\}>0\,.
\end{equation}
\end{pro}

\noindent
\proof
First we note that, since $\omega_+>\omega_->0$,
and we are assuming
$\omega_+\neq 3 \omega_-$, we get $\gamma>0$.
Since
$|\bal+\bbe|=4$, $\bal\neq\bbe$,
the quantity $|\bom\cdot(\bal-\bbe)|$ can be 
one of the following:
	$$
		4\omega_+,\  4\omega_-, \ 
		|3\omega_--\omega_+|,\  3\omega_+-\omega_-,\ 
		2(\omega_+-\omega_-), \ \omega_+-\omega_-,\ 
		2\omega_+,\  2\omega_-\,.
	$$
Since $\omega_+>\omega_->0$ we have that
$
4\omega_+, 4\omega_-, 
2\omega_+, 2\omega_->\omega_-\geq \gamma
$
and
$
3\omega_+-\omega_->
2(\omega_+-\omega_-)>\omega_+-\omega_-\geq \gamma;
$
finally, by definition, $|3\omega_--\omega_+|\geq \gamma$,
proving \eqref{stoccolma}.
\\
The resonant case and formula \eqref{stoccolma2} follow as well.
\eproof

\subsection{Resonant and nonresonant Birkhoff Normal Form}

As a consequence of the above Proposition \ref{dartagnan},
if $3\omega_-\neq\omega_+$,
we can solve the homological equation
\eqref{beige} since the coefficients in \eqref{vivaldi} are well defined.
If this is not the case, we can \emph{always} solve
the ``resonant homological equation'':
	\begin{equation}\label{beigeRES}
			\mathtt G+\{S,{\mathtt N}\}		=\bar{\mathtt H}_{4,{\rm res}}\,,
		\end{equation}	
		where
	\begin{equation}\label{H4RES}	
		\bar{\mathtt H}_{4,{\rm res}}
		:=
			\bar{\mathtt H}_4+
	\mathtt G_{(0,1),(3,0)}z_2\bar z_1^3+
	\mathtt G_{(3,0),(0,1)}z_1^3\bar z_2
	\stackrel{\eqref{Gab}}=
	\bar{\mathtt H}_4+
	\frac{f_{3,1}}
				{4(\sqrt{\omega_-})^{3}
				\sqrt{\omega_+}}
	(z_2 \bar z_1^3+z_1^3\bar z_2)
	\,.
	\end{equation}

Finally we note that
 the map  $\mathbf z\to\Phi^1_S(\mathbf z)$
 is defined in a neighborhood 
 of the origin and is close to the identity 
 since\footnote{Recalling Footnote \ref{footnote6}
the Hamiltonian flow $\Phi_S^t$ satisfies the equation
 $\frac{d}{dt}\Phi_S^t=
\mathrm i \mathbf J \nabla S(\Phi_S^t)$.
Since $S$ is a homogeneous polynomial of degree four, its
gradient has degree three, namely $|\nabla S|=O(\|\mathbf z\|^3)$.
}
$\Phi^1_S(\mathbf z)=\mathbf z+
O(\|\mathbf z\|^3)$, where
$$\|\mathbf  z\|=\|(z_1,z_2)\|:=\max\{|z_1|,|z_2|\}\,.
$$

\noindent
By denoting by $\Phi_{{\rm res}}=\Phi_S^1$
the time 1 Hamiltonian flow generated by   $S$ solving \eqref{beigeRES},
and by $\Phi_{\rm nonres}=\Phi_S^1$ the flow generated by
 $S$ solving \eqref{beige}
 we obtain  the following
	\begin{thm}[Birkhoff Normal Form]\label{rossoscuro}
			In a sufficiently small neighbourhood of the origin
			the symplectic transformation
			$\mathbf z\to\Phi_{{\rm res}}(\mathbf z)=\mathbf z+
			O(\|\mathbf z\|^3)$
			puts the Hamiltonian $\mathtt H$
			in resonant fourth order Birkhoff Normal Form, namely,
		\begin{equation}\label{giuliaRES}
			 \mathrm H_{{\rm res}}:=
			 \mathtt H\circ\Phi_{{\rm res}}
			 ={\mathtt N}+
			 \bar{\mathtt H}_{4,{\rm res}}+R_{\rm res}\,,\qquad
R_{\rm res}=O(\|\mathbf z\|^6)\,.
		\end{equation}
	 Moreover,
		if $\omega_+\neq 3 \omega_-$
		in a sufficiently small neighbourhood of the origin,
		the symplectic transformation
		$\mathbf z\to\Phi_{\rm nonres}
		(\mathbf z)=
				\mathbf z+O(\|\mathbf z\|^3)$
				yields the Hamiltonian $\mathtt H$
		in fourth order Birkhoff Normal Form, namely,
		\begin{equation}\label{giulia}
			 \mathrm H:=
			 \mathtt H\circ\Phi_{\rm nonres}
			 ={\mathtt N}+
			 \bar{\mathtt H}_{4}
			 +R\,,
			 \qquad		 
			 R=O(\|\mathbf z\|^6)\,.
		\end{equation}
 	\end{thm}

\subsection{Estimate on the remainder
of the BNF and on the applicability
threshold of the amplitudes}
\label{subsec:BNF}

We now  
estimate  the remainder
$R$
of the BNF in \eqref{giulia}
and  
 the applicability
threshold of the amplitudes,
see
\eqref{ampiezzeiniziali}, namely
the maximal values of the amplitudes
for which the BNF procedure 
makes sense.
We focus on the nonresonant case.
The corresponding estimates
for the resonant case are analogous
(even better, see Remark \ref{remy}).

\begin{pro}[Quantitative estimates
on nonresonant BNF]
\label{sopra}
Let $r>0$ and $0<\delta<1$.
Assume that $3 \omega_-\neq \omega_+$
and that $S$
in \eqref{beige} satisfies 
the smallness condition
\begin{equation}\label{tresca}
 S_r\leq (1-\delta)r\,, 
\qquad
S_r:=\max_{\|\mathbf  z\|\leq  r}
\|\nabla S\|\,.
\end{equation}
 Then the symplectic transformation
$\Phi_{\rm nonres}=\Phi_S^1$
is well defined
and close to the identity:
$$
\Phi_{\rm nonres}:B_{\delta r}\to  B_r\,,\qquad
\|\Phi_{\rm nonres}(\mathbf  z)-\mathbf  z\|\leq 
S_r\,.
$$
Furthermore
\begin{equation}\label{stimaresto}
\max_{\|\mathbf  z\|\leq \delta r}
|R|\leq 
\max_{\|\mathbf  z\|\leq r}
\Big(
|\{S,\bar{\mathtt H}_4\}|+\frac12
|\{S,\hat{\mathtt G}\}|
\Big)
\,,
\end{equation}
where $\hat{\mathtt G}:=\mathtt G-\bar{\mathtt H}_4$.
\end{pro}
\noindent
\proof
Denoting by $\mathbf z(t)$
the solution of Hamilton's  equations
of $S$ with initial datum 
$\mathbf z(0)=\mathbf z$,
one has that 
$\Phi_{\rm nonres}(\mathbf z)=
\Phi_S^1(\mathbf z)=\mathbf z(1)$.
By \eqref{tresca}, for every 
$0\leq t\leq 1$,
 the Hamiltonian flow 
$$
\Phi_S^t:
\{\|\mathbf  z\|\leq \delta r\}\to
\{\|\mathbf  z\|\leq r\}
$$
is well defined and
$$
\|\Phi_S^1(\mathbf  z)-\mathbf  z\|
=
\|\mathbf  z(1)-\mathbf  z(0)\|\leq 
S_r
$$
holds.
Finally, by \eqref{resto},
we get \eqref{stimaresto}.
\eproof

\medskip

We now rewrite Proposition \ref{sopra}
 in terms of the coefficients of
 $\mathtt G$ and $S$ namely,
  by \eqref{vivaldi},
 in terms of the coefficients of
 $\mathtt G$, $\omega_-$ and $\omega_+$.
First we need some notations.
 For a   polynomial $P$ of degree four
written in the form  \eqref{fuga}, we set
 \begin{equation}\label{25}
 P^{(j)}:=
 \sum_{|\bal+\bbe|=4}
  \alpha_j|P_{\bal,\bbe}|\,,\quad j=1,2,\qquad
  P_*:=\max\{P^{(1)},P^{(2)}\}\,.
 \end{equation}
 If $P$ has the property that
 $|P_{\bal,\bbe}|=|P_{\bbe,\bal}|$
 for every $\bal,\bbe$ with $|\bal+\bbe|=4$,
 then
  for $j=1,2$
\begin{equation}\label{flauto}
\max_{\|\mathbf  z\|\leq r}
|\partial_{z_j}P|
=
\max_{\|\mathbf  z\|\leq r}
|\partial_{\bar z_j}P|
\leq
r^3 P^{(j)}
\,.
\end{equation}
As a consequence, one can easily
estimate the gradient of a polynomial
on the complex domain 
$\{\|\mathbf  z\|\leq r\}$
 in terms of its coefficients:
 \begin{equation}\label{forchettata}
 \max_{\|\mathbf  z\|\leq r}
 \|\nabla P\|\leq P_* r^3\,.
 \end{equation}
If $P$ and $\tilde P$
are as above,
also
the Poisson brackets in 
\eqref{poisson}
are estimated 
in terms of the coefficients of the 
two polynomials:
\begin{equation}\label{linguine}
\max_{\|\mathbf  z\|\leq r}
|\{P,\tilde P\}|
\leq 
2 r^6
\sum_{j=1,2}
P^{(j)} \tilde P^{(j)}\,.
\end{equation}

\begin{cor}
Let $r>0$ and $0<\delta<1$. The smallness condition \eqref{tresca}
 is satisfied assuming the stronger condition
 \begin{equation}\label{cibilta}
r^2 S_*\leq 1-\delta\,,
\end{equation}
 where, by \eqref{25} and \eqref{vivaldi},
  \begin{equation}\label{S*}
   S_*:=\max\{S^{(1)},S^{(2)}\}\,,
   \qquad
 S^{(j)}:=
 \sum_{|\bal+\bbe|=4}
  \alpha_j\frac{|\mathtt G_{\bal,\bbe}|}
  {|\bom\cdot(\bal-\bbe)|}\,,\quad j=1,2
 \,.
 \end{equation}
 Moreover, one obtains the (weaker)
 estimates 
 \begin{equation}\label{maschera}
 \|\Phi(\mathbf  z)-\mathbf  z\|\leq 
r^3 S_*\,,\qquad\qquad
\max_{\|\mathbf  z\|\leq \delta r}
|R|\leq R_\dag
r^6\,,
\qquad
R_\dag:= \sum_{j=1,2}
S^{(j)}\Big(2\bar{\mathtt H}_4^{(j)}
+\hat{\mathtt G}^{(j)}
\Big)\,.
 \end{equation}
 \end{cor}

\begin{figure}[h!]
\center
\includegraphics[width=13cm,height=13cm,keepaspectratio]{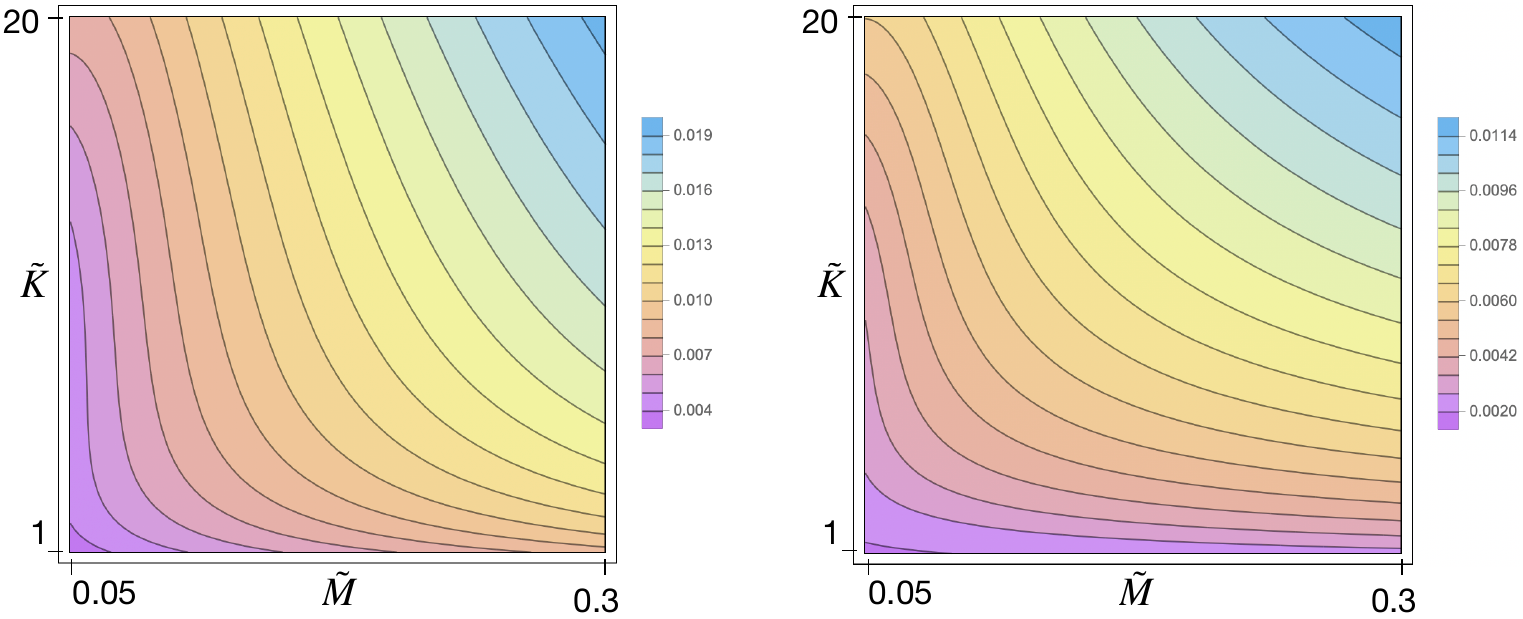}
\caption{Level curves of the admissible 
wave amplitude for the resonant BNF of the
optical (left) and acoustic (right) modes
as function of $\tilde M,\tilde K$
are plotted
in the case of the honeycomb metamaterial 
studied above.
Here $N_{3}=-10^4$
and the wave numbers are chosen
in $\mathbf X$. 
Compare with the analogous portraits in Figure 
\ref{Soglia_aPiuPar} for
the admissible amplitudes of the nonresonant
BNF.}
\label{Soglia_aPiuParRIS}
\end{figure}

	\begin{rem}\label{remy}
 Analogous estimates can be 
 obtained
 also in the resonant case.
 Since, in this case, the term 
 ${S}_{\bal,\bbe}$ with
 $(\bal,\bbe)=(3,0,0,1), (0,1,3,0)$
 does not appear in the definition of
 $S$ given in \eqref{vivaldi},
 the quantity $S_*$ is smaller
 and the applicability threshold
 for the amplitudes
 is larger; compare  Figures \ref{Soglia_aPiuPar} and  \ref{Soglia_aPiuParRIS}.
\end{rem}

\medskip
In terms of the initial datum
$\mathbf  z(0)$. Since $\|\mathbf  z(0)\|\leq \delta r$,
we can rewrite the smallness condition
\eqref{cibilta}
(ensuring
\eqref{tresca})
as 
$\|\mathbf  z(0)\|^2\leq (1-\delta)\delta^2/ S_*$.
One can optimise taking $\delta=2/3$,
which reaches the maximum of 
$(1-\delta)\delta^2$, namely $4/27$.
So the smallness condition
on the initial datum reads
\begin{equation}\label{soglia}
\epsilon:=\|\mathbf  z(0)\|\leq \frac{2}{3\sqrt{3S_*}}\,.
\end{equation}
Transforming  back into the original variables
$\mathbf q=(q_1,q_2)$
through
 \eqref{islanda} and \eqref{islanda2} yields
 \begin{equation}\label{esagerato}
q_1=(z_1+\bar z_1)/\sqrt{2\omega_-}\,,\qquad
q_2=(z_2+\bar z_2)/\sqrt{2\omega_+}\,,
\end{equation}
and the smallness condition
\eqref{soglia}
 on the
  initial amplitudes 
  $a_+$ and $a_-$
  reads:
\begin{equation}\label{ampiezzeiniziali}
|q_1(0)|\leq a_-=\sqrt{\frac{2}{27\omega_- S_*}}\,,
\qquad
|q_2(0)| \leq a_+=\sqrt{\frac{2}{27\omega_+ S_*}}\,.
\end{equation}
Let us rephrase 
condition \eqref{ampiezzeiniziali}
in terms of the detuning parameter
$\sigma=\omega_+ - 3 \omega_-$
and of the perturbative parameter
$\varepsilon$, which measures  the distance
from the origin in the original variables.
It is easy to choose a constant
$\tilde c>0$ such that
$|\sigma|\leq \tilde c \gamma$
(defined in \eqref{stoccolma}).
Then, by recalling 
 the definition of
$S_*$ given in \eqref{S*},
we see that
$S_*\leq c_\dag/|\sigma|$, where
$c_\dag$ depends only on $\tilde c$ and on the 
coefficients of $\mathtt G$
given in \eqref{calippo}.
Then we can define a constant
$c_1$ such that the smallness 
condition \eqref{ampiezzeiniziali}
reads, recalling \eqref{soglia} and \eqref{esagerato},
$$
\varepsilon\leq c_1\sqrt{|\sigma|}\,,
$$
justifying the formulas given in
the introduction (see page
\pageref{zoomarine}).

\subsection{Action-angle variables}

We now introduce action-angle variables\footnote{$\mathbb T:=\mathbb R/{2\pi \mathbb Z}$, $\mathbb T^2:=\mathbb R^2/{2\pi \mathbb Z^2}$.}
$(\mathbf I,\bphi)=(I_1,I_2,\varphi_1,\varphi_2)\in \mathbb R^2\times \mathbb T^2$
through the transformation
	\begin{equation}\label{amazon}
		 z_j=\sqrt {I_j} e^{-\mathrm i \varphi_j}\,,\qquad
			I_j>0\,,
			\qquad j=1,2.
	\end{equation}
	
\begin{rem}
 Note that the above map is singular at 
	$ z_1$ or $z_2=0$ 
and is defined for $I_1,I_2>0$.
\end{rem}
The change of variables in \eqref{amazon}
is symplectic:  a Hamiltonian
	$\mathtt H(\mathbf z,\bar{\mathbf z})$ with Hamilton's equations
	$\dot{\mathbf z}=-\mathrm i \partial_{\bar{\mathbf z}}\mathtt H$ and
	$\dot {\bar {\mathbf z}}=\mathrm i \partial_{\mathbf z}\mathtt H$,
	is transformed into a Hamiltonian
	$\mathcal H(\bI,\bphi):=\mathtt H(\sqrt \bI e^{-\mathrm i \bphi},
	\sqrt \bI e^{\mathrm i \bphi})$ with Hamilton's equations
$\dot \bI=-\partial_{\bphi}\mathcal H$ 
and $\dot\bphi=\partial_\bI \mathcal H$.
In particular $\bar{\mathtt H}_{4}$ and $\bar{\mathtt H}_{4,{\rm res}}$ in \eqref{mozart}, \eqref{H4RES} are transformed into\footnote{We use that if $a\in \mathbb C$\,, then $a+\bar a=2\mathrm{Re}(a)$ and $\theta \in \mathbb R$, $e^{\mathrm i \theta} =\cos \theta +\mathrm i \sin \theta$.}
	\begin{eqnarray}\label{pollo}
						\mathcal H_{4}(\bI)
					&:=&
						{\mathtt G}_{(2,0),(2,0)} 
						I_1^2
						+
						{\mathtt G}_{(1,1),(1,1)}
						I_1 I_2
							+
						{\mathtt G}_{(0,2),(0,2)} 
						I_2^2\,,
							\\
						\mathcal H_{4,{\rm res}}(\bI,\bphi)
					&:=&
						\mathcal H_{4}(\bI)+
						\frac{f_{3,1}}
						{2(\sqrt{\omega_-})^{3}
						\sqrt{\omega_+}}
						\sqrt{I_1^3I_2}
						\cos (\varphi_2-3\varphi_1)\,,
							\label{polloRES}
	\end{eqnarray}
respectively.
Recollecting the Hamiltonians 
in action angle variables
in the nonresonant, respectively,  resonant
case are:
	\begin{eqnarray}\label{HamAA}
					\mathcal H(\bI,\bphi)
					&=&\omega_-I_1+\omega_+ I_2
					+\mathcal H_{4}(\bI)+\mathcal R(\bI,\bphi)\,,
\qquad\ \ \ \ \qquad \mathcal R=O(\|\bI\|^3)\,,
							\\
					\mathcal H_{\rm res}(\bI,\bphi)
					&=&\omega_-I_1+\omega_+ I_2
					+\mathcal H_{4,{\rm res}}(\bI,\bphi)+
\mathcal R_{\rm res}(\bI,\bphi)\,,
\qquad \mathcal R_{\rm res}=O(\|\bI\|^3)					\,.
					\label{HamAARES}
	\end{eqnarray}
For small amplitude, namely for small values
of $\|\bI\|$ the term $O(\|\bI\|^3)$ is negligible
and we consider only the main term 
of the Hamiltonian:
\begin{eqnarray}\label{secular}
					\hat{\mathcal H}(\bI,\bphi)
					&=&\omega_-I_1+\omega_+ I_2
					+\mathcal H_{4}(\bI)\,,
					\nonumber
							\\
					\hat{\mathcal H}_{\rm res}(\bI,\bphi)
					&=&\omega_-I_1+\omega_+ I_2
					+\mathcal H_{4,{\rm res}}(\bI,\bphi)\,.
	\end{eqnarray}

\section{The nonresonant case}

In this section we consider  the nonresonant case
$\omega_+\neq 3\omega_-$.

\subsection{Nonlinear frequencies}

By \eqref{giulia} we are dealing with the Hamiltonian
$ \hat{\mathcal H}$ in \eqref{secular}
			 whose new frequencies
			 are
			 the derivatives
of the energy with respect to the actions, namely, by \eqref{pollo},
\begin{eqnarray}\label{omega_nonlineare}
			\omega_-^{\rm nl}
&:=&\partial_{I_1} \hat{\mathcal H}=	
\omega_-+
			2{\mathtt G}_{(2,0),(2,0)}I_1
			+
	{\mathtt G}_{(1,1),(1,1)}I_2
	\nonumber
	\\
	&=&\omega_-+
	\frac{3}{4\omega_-^2}
				\Big(
				(\phi_1^-)^4M_3+(\phi_2^-)^4N_3
				\Big)I_1
				\nonumber
				\\
				&&\quad\ \,
				+
				\frac3{2\omega_- \omega_+}\Big(
				(\phi_1^-)^2(\phi_1^+)^2M_3+(\phi_2^-)^2(\phi_2^+)^2 N_3
				\Big)I_2
				\,,
				\nonumber
					\\
			\omega_+^{\rm nl}
			&:=&\partial_{I_2} \hat{\mathcal H}
			=
			\omega_++2{\mathtt G}_{(0,2),(0,2)}I_2+
{\mathtt G}_{(1,1),(1,1)}I_1
\nonumber
	\\
	&=&\omega_+ +
	\frac{3}{4\omega_+^2}
				\Big(
				(\phi_1^+)^4M_3+(\phi_2^+)^4N_3
				\Big)I_2
				\nonumber
				\\
				&&\quad\ \,
				+
				\frac3{2\omega_- \omega_+}\Big(
				(\phi_1^-)^2(\phi_1^+)^2M_3+(\phi_2^-)^2(\phi_2^+)^2 N_3
				\Big)I_1
				\,,
	\end{eqnarray}
by \eqref{calippo}.	
In particular, when 
$M_3=0$, we 
have
\begin{eqnarray}\label{omega_nonlineareN3}
			\omega_-^{\rm nl}
&=&		
\omega_-+N_3\left(
	\frac{3}{4\omega_-^2}
				(\phi_2^-)^4 I_1
				+
				\frac3{2\omega_- \omega_+}
			(\phi_2^-)^2(\phi_2^+)^2
				I_2
				\right)
				\,,
				\nonumber
					\\
			\omega_+^{\rm nl}
			&=&
			\omega_+ +N_3\left(
	\frac{3}{4\omega_+^2}
				(\phi_2^+)^4
				I_2
				+
				\frac3{2\omega_- \omega_+}(\phi_2^-)^2(\phi_2^+)^2 I_1
				\right)
				\,.
	\end{eqnarray}
By \eqref{islanda} and \eqref{islanda2}
we have that the amplitudes $a_-, a_+$
in the original variables $q_1$ and $q_2$,
namely when 
\begin{equation}\label{cerniera}
q_1(0)=a_-\,,\quad 
q_2(0)=a_+\,,\quad
\dot p_1(0)=\dot p_2(0)=0\,,
\end{equation}
correspond, by \eqref{amazon}, to initial 
action-angle variables: 
\begin{equation}\label{ampiezze}
I_1(0)=\frac12\omega_- a_-^2\,,\qquad
I_2(0)=\frac12\omega_+ a_+^2\,,\qquad
\varphi_1(0)=\varphi_2(0)=0\,.
\end{equation}
Therefore  \eqref{omega_nonlineareN3} is exactly
\eqref{formulaSL1}.

\begin{rem}\label{remo}
Formula \eqref{omega_nonlineareN3},
equivalently \eqref{formulaSL1},
does not hold true close to resonances.
Obtaining the counterpart of formula  
\eqref{omega_nonlineareN3} in the resonant case is much
more complicated since one has to 
integrate the Hamiltonian
$\hat{\mathcal H}_{\rm res}$ in \eqref{secular}.
This approach involves introducing appropriate action-angle coordinates using elliptic integrals. While the topic is well-established, deriving explicit analytic formulas that account for all the physical parameters is challenging. A more detailed explanation will be provided in a forthcoming work
\cite{DL2}.
\end{rem}

	Note that the Hamiltonian in  the first line of
	\eqref{secular}
	is \emph{integrable}
	since the actions $I_1$ and $I_2$
	are constants of motion.
	The solutions of
	 Hamilton's equations are
	\begin{equation*}\label{HEAA}
			I_1(t)=I_1(0)\,,
			\ \ \
			I_2(t)=I_2(0)\,,
			\ \ \
			\varphi_1(t)= \varphi_1(0)+
			\omega_-^{\textrm{nl}} t\,,
			\ \ \
			\varphi_2(t)= \varphi_2(0)+
			\omega_+^{\textrm{nl}} t\,.
	\end{equation*}

\begin{rem}\label{rem:KAM}
The truncated Hamiltonian  $\hat{\mathcal H}$
in  the first line of \eqref{secular} is integrable and
all its solutions are quasi--periodic in time with 
frequencies $\omega_\pm^{\rm nl}$
given by \eqref{formulaSL1}.
This is not the case of the complete Hamiltonian
 $\mathcal H$ in \eqref{HamAA}.
Since the action-to-frequency map
 $(I_1,I_2)\to (\omega_-^{\rm nl},\omega_+^{\rm nl})$ in \eqref{omega_nonlineareN3} is invertible
 (its Jacobian determinant being 
 $-(9/4)\omega_-^{-2} \omega_+^{-2}
			(\phi_2^-)^4(\phi_2^+)^4<0$),
the KAM Theorem would ensure that
the majority (in the sense of measure)
of solutions of $\mathcal H$ with small  initial actions
(equivalently, small initial   amplitudes)
are still quasi-periodic
with (strongly irrational) frequencies close to
$\omega_\pm$. 
Every quasi-periodic orbit densely covers
 the surface of a two dimensional 
 torus, which is  invariant for the Hamiltonian
 flow of $\mathcal H$, and is stable,
 in the sense that the values of the actions, as well as of the amplitudes, remain close to the
 initial ones for all times.
 As well known,
 a remarkable byproduct of the KAM Theorem\footnote{When the system,
 depending on the parameters, satisfies the hypotheses of the KAM Theorem.}
 in two degrees of freedom is that \textsl{all}
 solutions, even those that are not quasi-periodic, are perpetually stable.
 Indeed, by energy conservation, 
 every orbit is confined to its 
 three-dimensional energy level surface.
 The above two-dimensional invariant tori separate the energy level, and every orbit either lies on an invariant torus or is trapped among two of them. In both cases no escape is possible, and the actions stay forever close to their initial values.
\end{rem}

\subsection{Asymptotic solutions 
and estimates on the remainder}

Here, following \cite{SW23jsv}, we consider the 
multiple scales
development   of the solution
in Cartesian variables up to the third
order in the perturbative parameter,
see \eqref{tommaso}
and \eqref{sviluppo},
 and give explicit
estimates on the 
fifth order
remainder, see \eqref{stivali}.

Let us consider the solutions of Hamilton's  
equation 
$\dot {\mathbf z}=-\rm i \partial_{\mathbf z}
\mathrm H$
of the 
Hamiltonian
$
			 \mathrm H
			 ={\mathtt N}+
			 \bar{\mathtt H}_{4}
			 +R
			 $
in \eqref{giulia}.
We assume that the initial datum
$\mathbf z(0)$ is
close to the origin, namely the parameter $\epsilon=\|\mathbf z(0)\|$
introduced in \eqref{soglia} is small.
\\
Let $0<\xi\leq 1/2$ such that 
\begin{equation}\label{parigi}
2\xi\epsilon=\min\{|z_1(0)|,|z_2(0)|\}\,.
\end{equation}
As done in the multiple scales approach we
  write the solution in the form\begin{equation}\label{tommaso}
\mathbf z(t)=\epsilon \mathbf z^{(1)}(t,\epsilon^2 t)
+\epsilon^3 \mathbf z^{(3)}(t,\epsilon^2 t)
+\epsilon^5\tilde{\mathbf z}(t)\,,
\end{equation}
 where $\mathbf z^{(1)}(T_0,T_2)$ and
$\mathbf z^{(3)}(T_0,T_2)$ are periodic
(with possibly different periods) in both
their arguments; moreover, 
we give explicit  estimates on 
 $\|\tilde{\mathbf z}(t)\|$
 in terms of the coefficients of
 $\mathtt G$, $\omega_-$ and $\omega_+$.
 \\
 Let us consider the complex domain
 $\mathcal D:=\{\xi\epsilon\leq |z_1|\leq 2\epsilon\}\times
 \{\xi\epsilon\leq |z_2|\leq 2\epsilon\}$,
 which is the product of two annuli.
Recalling the definition of $R_\dag$, $S_*$ and
$\mathcal H_{4}$
given in, respectively, \eqref{maschera}, \eqref{S*} and \eqref{pollo},
we state the following
\begin{pro}\label{caffe}
Let 
\begin{equation}\label{osiride}
R_*:=(9/2)^6 R_\dag\,, \quad 
t_*:=\frac{3\xi^2}{2R_*}\,,
\quad
\mbox{and}\quad
\mathbf A:=\partial^2_\bI \mathcal H_{4}(0)\,,
\end{equation}
namely the Hessian matrix of 
$\mathcal H_{4}(\bI)$ evaluated at $\bI=0$.
 Assume the smallness condition
\begin{equation}\label{notredame}
\epsilon\leq \frac{2}{9\sqrt{3S_*}} \,.
\end{equation}
Then
$\mathbf z^{(1)}(T_0,T_2)=
\epsilon^{-1}\mathbf z(0)
e^{-\rm i (\bom T_0+\mathbf A\bI_0T_2)}$ and $\mathbf z^{(3)}\equiv 0$ and
\begin{equation}\label{stivali}
\|\mathbf z(t)\|\leq
 \frac{7}{2\sqrt 2} \, \epsilon
\,,\qquad
\|\tilde{\mathbf z}(t)\|
\leq
\left(\epsilon^2\|\mathbf A\||t|+
\xi^{-1}
\right)R_*|t|
\,,\qquad
\forall\, 
|t|\leq t_*
\epsilon^{-4}
\,.
\end{equation}
\end{pro}

\begin{rem}
 The second estimate in \eqref{stivali}
 is relevant only for $|t|\ll \epsilon^{-3}$
 so that 
 $\epsilon^5\|\tilde{\mathbf z}(t)\|\ll\epsilon$.
 Anyway
 $\epsilon^5\|\tilde{\mathbf z}(t)\|\leq \|\mathbf z(t)-\epsilon \mathbf z^{(1)} \|
 \leq \|\mathbf z(t)\|+\epsilon\leq  (1+7/2\sqrt 2) \, \epsilon$,
for every $
|t|\leq t_*
\epsilon^{-4}$.
\end{rem}

The proof of Proposition \ref{caffe}
can be found in the Appendix.


\section{Conclusions}


This study presented an in-depth analysis of the wave propagation equations governing 2D metamaterials with embedded periodically distributed Duffing-type resonators.
The system  can be generalized into a system of two coupled harmonic oscillators with cubic nonlinearity. The  system is Hamiltonian, with an elliptic equilibrium at the origin and two distinct linear frequencies associated with the acoustic and optical modes. 

In the small amplitude regime, away from resonances, Hamiltonian Perturbation Theory was employed to derive the first-order nonlinear correction to the linear frequencies. We also analytically estimated the threshold for the perturbative approach, defining the maximum amplitude at which these corrections remain valid. This threshold is critical for determining the range of applicability of the obtained analytical results.

Our findings show that this threshold decreases significantly in the presence of resonances, particularly when the ratio between the optical and acoustic frequencies approaches 3. Notably, the three-to-one internal resonance is the only first-order resonance that significantly affects the nonlinear frequency correction, highlighting its strong influence on the system's dynamics. To refine our understanding, we also evaluated the remainder after the first step of the perturbation procedure, which is essential for assessing the accuracy of the perturbation theory.

As a practical application of our theoretical results, we analyzed the wave propagation problem and the emergence of bandgaps  in metamaterial honeycombs with periodically distributed nonlinear resonators. This analysis provided valuable quantitative insights into how internal resonances impact the nonlinear appearance of  the bandgap and whether the width of the bandgap gets enhanced or demoted as a result of the nonlinear local resonances. Specifically, we identified a 3:1 resonant curve in the $(\tilde M,\tilde K)$-plane. This curve divides the parameter space into two regions: one where the dynamics are primarily nonresonant, allowing for the use of known formulas for nonlinear frequencies, and another where resonance effects dominate, requiring new expressions for the resonant nonlinear frequencies. In the nonresonant region, we identified a large area where the nonlinear bandgap, in the softening case, is significantly larger than the linear one.


%
%
%
%
%

\section{Appendix}

\subsubsection*{Proof of Proposition \ref{caffe}}

Taking $\delta=2/3$ and 
$r=3\epsilon/\delta=9\epsilon/2$ in \eqref{cibilta},
the smallness condition \eqref{cibilta}
is equivalent to \eqref{notredame}.
By \eqref{maschera} we have
\begin{equation}\label{cartagine}
\max_{\|\mathbf  z\|\leq 3\epsilon}
|R|\leq R_*
\epsilon^6
\,.
\end{equation}
Since $R(\mathbf z)$ is holomorphic
on the complex ball $\{\|\mathbf  z\|\leq 3\epsilon\}$,
by Cauchy estimates
we obtain
\begin{equation}\label{lenticchia}
\max_{\mathcal D}|\partial_{z_j} R(\mathbf z)|
\leq R_*
\epsilon^5\,,\qquad j=1,2\,.
\end{equation}
Now it is convenient to 
write $\mathrm H$ in action angle variables
as in 
\eqref{HamAA}, namely
$$
\mathcal H(\bI,\bphi)
=\bom \cdot \bI
+\frac12\mathbf A \bI\cdot\bI
+\mathcal R(\bI,\bphi)
$$	
where $\mathbf A$ is the Hessian matrix
of $\mathcal H_{4}(\bI)$ evaluated at $\bI=\mathbf 0$.	
Since by \eqref{amazon}
$I_j=|z_j|^2$, $j=1,2$,
the domain $\mathcal D$ in the action
reads
$\tilde{\mathcal D}:=
\{\xi^2\epsilon^2\leq I_1\leq 4\epsilon^2\}\times
 \{\xi^2\epsilon^2\leq I_2\leq 4\epsilon^2\}$.
By recalling \eqref{amazon},
 we obtain, by definition, 
$\mathcal R(\bI,\bphi)=
R(\sqrt I_1 e^{-\mathrm i \varphi_1},\sqrt I_2 e^{-\mathrm i \varphi_2})$.
Then, by 	\eqref{lenticchia},
we get for $j=1,2$,
\begin{equation}\label{lenticchia2}
\max_{\|\bI\|\leq 4\epsilon^2}|\partial_{\varphi_j} 
\mathcal R(\bI,\bphi)|
\leq 2R_*
\epsilon^6\,,\qquad 
\max_{\tilde{\mathcal D}}|\partial_{I_j} 
\mathcal R(\bI,\bphi)|
\leq \frac{R_*}{2\xi}
\epsilon^4
\,.
\end{equation}	 
Hamilton's  equations are	
expressed as	 
\begin{equation}\label{plumcake}
\dot{\bI}=-\partial_{\bphi}\mathcal H(\bI,\bphi)
=-\partial_{\bphi}\mathcal R(\bI,\bphi)
\,,
\qquad
\dot{\bphi}=\partial_{\bI}\mathcal H(\bI,\bphi)
=\bom+\mathbf A\bI+\partial_{\bI}\mathcal R(\bI,\bphi)
\end{equation}		 
with initial data $\bI(0)=:\epsilon^2 \bI_0$
(recall  \eqref{amazon} and  \eqref{soglia}),
$\bphi(0)=\bphi_0$, such that
$\mathbf z(0)=\epsilon\sqrt{\bI_0}e^{-\rm i\bphi_0}$.	
Note that by definition
\begin{equation}\label{prugna}
\|\bI_0\|=1\,,\qquad \min\{|I_{01}|,|I_{02}|\}=4\xi^2\leq 1\,.
\end{equation}
In view of \eqref{lenticchia2} we look for solutions $\bI(t),\bphi(t)$
of the form
\begin{equation}\label{splendente}
\bI(t)=\epsilon^2\bI_0+\epsilon^6 \bJ(t)\,,
\qquad
\bphi(t)=\bphi_0
+(\bom+\epsilon^2\mathbf A\bI_0)t
+\epsilon^4\bpsi(t)
\,.
\end{equation}
Note that $t\to \big(\epsilon^2\bI_0,\bphi_0
+(\bom+\epsilon^2\mathbf A\bI_0)t\big)$
is the solution of the truncated Hamiltonian
$\bom \cdot \bI+
\mathcal H_{4}(\bI)=
\bom \cdot \bI
+\frac12\mathbf A \bI\cdot\bI$
with initial datum	
$(\epsilon^2 \bI_0,\bphi_0)$.
By \eqref{plumcake} and \eqref{splendente}
 the higher order correction
$\big(\bJ(t),\bpsi(t)\big)$ satisfies
\begin{equation}\label{sparta}
\dot{\bJ}(t)
=-\epsilon^{-6}\partial_{\bphi}\mathcal R
\big(\bI(t),\bphi(t)\big)
\,,
\qquad
\dot{\bpsi}
=\epsilon^2\mathbf A\bJ+
\epsilon^{-4}\partial_{\bI}\mathcal R
\big(\bI(t),\bphi(t)\big)\,,
\end{equation}
with initial datum $\bJ(0)=\bpsi(0)=\mathbf 0$.
We claim that
\begin{equation}\label{atene}
|t|\leq t_*\epsilon^{-4}
\qquad
\Longrightarrow
\qquad
|\bJ(t)|\leq 2 R_* |t|\,,
\quad
|\bpsi(t)|\leq 
\left(\epsilon^2\|\mathbf A\| |t|+\frac1{2\xi}
\right)R_*|t|\,.
\end{equation}
In order to prove \eqref{atene}
we first show that, as long as
$|t|\leq t_*\epsilon^{-4}$, we have 
$\bI(t)\in \tilde{\mathcal D}$.
Indeed by contradiction  if 
there exists $t_0$ with $|t_0|<t_*\epsilon^{-4}$
such that
 $\bI(t)\in \tilde{\mathcal D}$
 for every $|t|\leq t_0$
 and $\bI(t_0)\in \partial\tilde{\mathcal D}$,
 namely
  $\|\bI(t_0)-\epsilon^2\bI_0\|\geq 3 \xi^2 \epsilon^2$ (recall \eqref{prugna}).
  Then, by \eqref{sparta}, \eqref{lenticchia2} and \eqref{osiride},
 we obtain
 $\|\bJ(t_0)\|\leq 2 R_* |t_0|<2R_*t_*\epsilon^{-4}=
 3\xi^2\epsilon^{-4}$.
Then by \eqref{splendente}
$$
3 \xi^2 \epsilon^2\leq
\|\bI(t_0)-\epsilon^2\bI_0\|=
\epsilon^6\|\bJ(t_0)\|<
 3\xi^2\epsilon^2\,,
$$
which is a contradiction.
In conclusion, we have proved that
$\bI(t)\in \tilde{\mathcal D}$
as long as
$|t|\leq t_*\epsilon^{-4}$.
Finally the claim in \eqref{atene}
follows by \eqref{sparta} and
\eqref{lenticchia2}.

Transforming  back to the variable 
$\mathbf z$, by \eqref{amazon}
and \eqref{splendente}
 we have that $\mathbf z(t)$
 in \eqref{tommaso}
is expressed as 
\begin{equation}\label{insetto}
\mathbf z(t)
=\sqrt{\bI(t)}e^{-\rm i \bphi(t)}
=
\epsilon
\sqrt{\bI_0+\epsilon^4 \bJ(t)}
e^{-\rm i (\bphi_0
+(\bom+\epsilon^2\mathbf A\bI_0)t)}
e^{-\rm i \epsilon^4\bpsi(t)}\,.
\end{equation}
Then we get 
\begin{equation}\label{sviluppo}
\mathbf z^{(1)}(T_0,T_2)=
\sqrt{\bI_0}
e^{-\rm i (\bphi_0
+\bom T_0+\mathbf A\bI_0T_2)}\,,\qquad
\mathbf z^{(3)}(T_0,T_2)=
\mathbf 0\,.
\end{equation}
By \eqref{insetto}, \eqref{prugna}, \eqref{atene}
and \eqref{osiride}
we get
$$
\|\mathbf z\|\leq \epsilon \sqrt 2(1+2R_* t_*)
\leq
\epsilon 7/2\sqrt 2 \,,\qquad
\forall\, |t|\leq t_*\epsilon^{-4}\,,
$$
proving the first estimate in \eqref{stivali}.
Regarding the fifth order remainder 
in \eqref{tommaso},
 for $j=1,2$ and $|t|\leq t_*\epsilon^{-4}$,
 by \eqref{insetto}
 we have\footnote{Using that
$|e^{\mathrm i x}-1|\leq |x|$ for $x\in\mathbb R$ and $|\sqrt{1+x}-1|\leq |x|/2$ for 
$x\geq -1$.}
\begin{eqnarray*}
&&\big|\tilde{z}_j(t)\big|
=\epsilon^{-5}
|z_j(t)-\epsilon z^{(1)}_j(t,\epsilon^2 t)|
\leq 
\epsilon^{-4}
\Big|\sqrt{I_{0j}+\epsilon^4 J_j(t)}
e^{-\rm i \epsilon^4\psi_j(t)}
-\sqrt{I_{0j}}
\Big|
\\
&\leq&
\epsilon^{-4}
\Big|\Big(\sqrt{I_{0j}+\epsilon^4 J_j(t)}
-\sqrt{I_{0j}}\Big)
e^{-\rm i \epsilon^4\psi_j(t)}
+\sqrt{I_{0j}}(e^{-\rm i \epsilon^4\psi_j(t)}-1)
\Big|
\\
&\leq&
\epsilon^{-4}
\Big|\sqrt{I_{0j}+\epsilon^4 J_j(t)}
-\sqrt{I_{0j}}
\Big|
+
\epsilon^{-4}
\sqrt{I_{0j}}
\big|e^{-\rm i \epsilon^4\psi_j(t)}-1\big|
\\
&\leq&
\frac{|J_j(t)|}{2\sqrt{I_{0j}}}
+\sqrt{I_{0j}}|\psi_j(t)|
\stackrel{\eqref{prugna}}\leq
\frac{1}{4\xi}|J_j(t)|+|\psi_j(t)|\,.
\end{eqnarray*}
Finally, by  \eqref{atene}, 
we obtain the second estimate in \eqref{stivali}.
\eproof

\end{document}